\documentclass[reqno, 10pt]{amsart}
\usepackage{amssymb, amsmath, amsthm, epsf, epsfig, xcolor}
\usepackage{verbatim}
\usepackage[shortlabels]{enumitem}
\usepackage{cite}

\newtheorem{theorem}{Theorem}
\newtheorem{proposition}{Proposition}%
\newtheorem{lemma}{Lemma}%
\newtheorem{conjecture}{Conjecture}%
\newtheorem{example}{Example}%
\newtheorem{remark}{Remark}%
\newtheorem{definition}{Definition}%

\newcommand{\beq}{\begin{equation}} 
\newcommand{\eeq}{\end{equation}} 
\newcommand{\bal}{\begin{align}} 
\newcommand{\eal}{\end{align}} 
\newcommand{\bals}{\begin{align*}} 
\newcommand{\eals}{\end{align*}} 
\newcommand{\overliner}[1]{\begin{array}{#1}} 
\newcommand{\earr}{\end{array}}

\newcommand{\bth}{\begin{theo}} 
\newcommand{\bl}{\begin{lemma}} 
\newcommand{\el}{\end{lemma}} 
\newcommand{\bp}{\begin{prop}} 
\newcommand{\ep}{\end{prop}} 
\newcommand{\bdf}{\begin{df}} 
\newcommand{\edf}{\end{df}} 
\newcommand{\brem}{\begin{rem}} 
\newcommand{\erem}{\end{rem}} 
\newcommand{\bnrem}{\begin{nrem}} 
\newcommand{\enrem}{\end{nrem}} 
\newcommand{\bex}{\begin{ex}} 
\newcommand{\eex}{\end{ex}} 
\newcommand{\bcor}{\begin{cor}} 
\newcommand{\ecor}{\end{cor}} 
\newcommand{\bncor}{\begin{ncor}} 
\newcommand{\encor}{\end{ncor}} 
\newcommand{\bpf}{\begin{proof}} 
\newcommand{\epf}{\end{proof}}

\def\({\left(} 
\def\){\right)}

\newcommand{\diff}{\Delta}

\numberwithin{equation}{section} 

\begin{document}
\title[Universal asymptotics of catalytic equations]{Universal asymptotic properties of positive functional equations with one catalytic variable}
\author{Michael Drmota$^*$, Eva-Maria Hainzl$^*$}

\thanks{{}$^*$ TU Wien, Institute of Discrete Mathematics and Geometry,
Wiedner Hauptstrasse 8-10, A-1040 Vienna, Austria. michael.drmota@tuwien.ac.at. Research 
supported by the
Austrian Science Foundation FWF, projects F50-02, F55-02, and P35016.}

\title[Universal asymptotics of catalytic equations]{Universal asymptotic properties of positive functional equations with one catalytic variable}

\begin{abstract}Functional equations with one catalytic appear in several combinatorial applications, most notably in the 
enumeration of lattice paths and in the enumeration  of planar maps. The main purpose of this 
paper is to show that under certain positivity assumptions the dominant singularity of the solutions 
function have a universal behavior. 
We have to distinguish between linear catalytic equations, where a dominating square root singularity
appears,  and non-linear catalytic equations, where we -- usually -- have a singularity of type $3/2$.
\end{abstract}

\maketitle

\section{Introduction}
\label{sec1}

Functional equations with one catalytic variable have their origin mostly in map enumeration \cite{Tutte1963} and 
in lattice path enumeration \cite{BandFla2002}. Such equations were first solved with
the help of the {\it kernel method} \cite{BandFla2002,Prodinger2003,AsinowskiBBG2020} (in the linear case) and with the help of the
{\it quadratic method} \cite{Tutte1963,BrownTutte} (in the quadratic case). Both approaches were unified and extended
by Bousquet-M\'elou and Jehanne \cite{BM-J}. They considered {\it general catalytic equations} of 
the form 
\begin{equation}\label{eqcatBMJ}
P(z,u, M(z,u), M_0(z), \ldots, M_{k-1}(z)) = 0,
\end{equation}
where $P(z,u,x_{0},x_{1},\dots x_{k})$ is a polynomial and all power series $M(z,u), M_0(z), \ldots, M_{k-1}(z)$ are characterized by this equation.

The variable ``$u$'' is called {\it catalytic} since it is usually an auxiliary variable that counts an additional 
(usually combinatorial) parameter which simplifies the recursive decomposition of the structure of interest. In general, one is just interested in the function $M(z,0)$, $M(z,1)$ or in $M_0(z)$.

One of the most prominent examples is the counting problem of rooted planar maps that goes back to Tutte \cite{Tutte1963}.
Let $M_k(z)$ denote the generating function of those maps, where the root face has valency $k\ge 0$.
Then we have $M_0(z) = 1$ and
\begin{align}
M_k(z) &= z \sum_{j=0}^{k-2} M_j(z)M_{k-j-2}(z) + z \sum_{j=k-1}^\infty M_j(z) \qquad (k\ge 1)  \label{eqmapsrec}
\end{align} 
{where the right sum arises if the root edge is not a bridge and the left sum if the deletion of the root edge results in decomposing the map into two components.}
One is interested in the generating function $M(z) = \sum_{k\ge 0} M_k(z)$ of all planar maps.
By introducing the variable~$u$ and setting  $M(z,u) = \sum_{k\ge 0} M_k(z) u^k$, the infinite system (\ref{eqmapsrec})
rewrites to the {\it catalytic equation} 
\begin{equation}\label{eqMzv}
M(z,u) = 1 + zu^2 M(z,u)^2 + uz \frac{ u M(z,u) - M(z,1) }{u-1}.
\end{equation}
By using the above mentioned quadratic method the equation can be explicitly solved:
\[
M(z) = M(z,1) = \frac{18z-1 + (1-12 z)^{3/2}} {54 z^2}.
\]
This leads to an explicit formula $M_n = [z^n]\, M(z,1) = \frac{2(2n)!}{(n+2)! n!} 3^n$ 
and to an asymptotic one: $M_n \sim (2/\sqrt \pi) 12^n n^{-5/2}$.
Note that the asymptotic behavior is reflected by the dominant singular behavior of 
$M(z,1)$ at $z_0 = 1/12$. 
We say that an algebraic function $f$ has a square root singularity at $z_0$ or that the singularity is of type $1/2$ if the dominating term in the
Puiseux expansion at $z_0$ is of the form $(z-z_0)^{1/2}$. That is,
\[
f(z) = f_0 + f_1 (z-z_0)^{1/2} + f_2(z-z_0) + f_3 (z-z_0)^{3/2}+ \cdots,
\]
where $f_1\ne 0$. Conversely if $f_1 = 0$ and $f_3\ne 0$,
then we say that the singularity at $z_0$ has type $3/2$.
{In the case of planar maps, the type of the singularity is indeed $3/2$ which translates to the critical exponent $-5/2 = -1-3/2$ by the well known Transfer Lemma 
\cite{FO}.}

In \cite{BM-J} several applications mostly from map enumeration (different classes
of planar maps, constellations, hard particles in planar maps etc.) are given.
Bousquet-M\'elou and Jehanne~\cite{BM-J} considered in particular equations of the form\footnote
{In \cite{BM-J} an additional additive polynomial term $F_0(z,u)$ is considered, however, by
substituting $M(z,u)$ by $M(z,u)-F_0(z,u)$ the seemingly more general case can be 
reduced to the case (\ref{eqBMJgen}).}
\begin{equation}\label{eqBMJgen}
M(z,u) = z Q\left(z,u, M(z,u), \diff^{(1)}(M(z,u)) , \ldots, \diff^{(k)}(M(z,u)) \right),
\end{equation}
where  $Q(z,u,\alpha_{0},\alpha_{1},\dots, \alpha_{k})$ is a polynomial and where we have used the notation $M(z,u) = \sum_{k\geq 0} M_k(z)u^k$ and 
\[
\diff^{(j)}(M(z,u)) =  \frac{ M(z,u) - M_0(z) - u\, M_{1}(z) - \cdots u^{j-1} M_{{j-1}}(z) }{u^j} 
\]
It is convenient to consider just the catalytic variable $u$ at $0$. In the above case of planar maps we 
substitute $u$ by $u+1$ to reduce it to this case.

One main result of \cite{BM-J} is that equations of the form (\ref{eqcatBMJ}) can be solved 
with the help of  proper systems of polynomial equations. Hence, the solutions are always algebraic functions
and consequently for every singularity we have a Puiseux expansion. However, this approach does not specify 
the kind of the Puiseux expansion. There is in principle no restriction on the rational exponents that might occur. 

However, if we consider the special case $k=1$ (in (\ref{eqBMJgen})) 
\begin{equation}\label{eqBMJ}
M(z,u) = z Q\left(z,u, M(z,u), \frac{ M(z,u)-M_0(z)}u \right),
\end{equation}
where $Q(z,u,\alpha_0,\alpha_1)$ is a polynomial with {\it non-negative coefficients},
then Drmota, Noy, and Yu \cite{DNY2021} showed that
there is a dichotomy (under natural conditions on $Q$).
If $Q$ is linear in $\alpha_0$ and $\alpha_1$ then the dominant singularity is of type $\frac 12$, that is, 
a square root singularity which leads to an asymptotic behavior for the coefficients of the form $\sim c \rho^n n^{-3/2}$.
However, in the non-linear case the dominant singularity is of type $\frac 32$ (as in the above mentioned example of
planar maps)  which means that the coefficients are asymptotically of the form $\sim c \rho^n n^{-5/2}$. 

In what follows we will focus on the case $k=2$ and show that the results for the case $k=1$ can be extended. However, there are several (major) differences. 
Whereas in the case $k=1$ the catalytic equation can be solved with the help of
a so-called positive system of polynomial equations (see \cite{BD,DNY2021}) which determines directly a dominant square root singularity
for the involved solution functions. This property is widely lost for the cases $k\ge 2$. Thus, it is necessary
to develop new methods and concepts in order to deduce the universal singular behavior. 
Clearly we expect similar properties for all equations of the form~(\ref{eqBMJgen}) as well as for 
systems of positive catalytic equations but the cases $k > 2$ are even more involved.


\section{Main results} \label{sec:generic}

Before we explain the solution method by Bousquet-M\'elou and Jehanne \cite{BM-J} for an equation of the form~(\ref{eqcatBMJ}), let us introduce the classic short-hand notation $f_z$ to denote the partial derivative $\frac{\partial f}{\partial z}$ of a function $f$ with respect to $z$. \\
Now recall that all unknown functions $M(z,u)$, $M_i(z)$ are uniquely defined by the equation. Thus we may assume they are already given and we consider the partial derivative with respect to $u$
\begin{align}\label{eq:partial_eq}
\frac{\partial}{\partial u}&P(z,u,M(z,u),M_0(z),\dots,M_{k-1}(z)) \\
&= P_u(z,u,M(z,u),M_0(z),\dots,M_{k-1}(z)) \\
&\qquad + P_{x_0}(z,u,M(z,u),M_0(z),\dots,M_{k-1}(z))\cdot M_u(z,u).
\end{align}
The equation 
\[P_{x_0}(z,u,M(z,u),M_0(z),\dots,M_{k-1}(z)) = 0\]
has exactly $k$ Puiseux series solutions $u_i(z), i=1,\dots k$ with $u_i(0) = 0$, which we assume to be distinct (otherwise one can perturb the system). As equation~(\ref{eq:partial_eq}) always has to be satisfied, $P_u = 0$ along these solutions and we obtain the algebraic system of $3k$ equations
\begin{align}
P(z,u_i(z), f_i(z), M_0(z), \ldots, M_{k-1}(z)) & = 0, \quad 1\le i \le k,   \nonumber \\
P_{x_0}(z,u_i(z), f_i(z), M_0(z), \ldots, M_{k-1}(z)) & = 0, \quad 1\le i \le k,  \label{eqBMJsystem} \\
P_u(z,u_i(z), f_i(z), M_0(z), \ldots, M_{k-1}(z)) & = 0, \quad 1\le i \le k, \nonumber
\end{align}
for the $3k$ unknown functions $u_1(z), \ldots, u_k(z)$, $f_1(z), \ldots, f_k(z)$, $M_0(z), \ldots, M_{k-1}(z)$ which take the value $0$ at the origin, where $f_i(z) = M(z,u_i(z))$ for all $1\leq i \leq k$ and the solutions lead to the unknown functions.\\

Before we consider the case $k=2$ in (\ref{eqBMJgen}) we comment on the case $k=1$.
(We recall that this case was fully discussed in \cite{DNY2021}). It is convenient to replace $M(z,u)$ by 
$M(z,u) =  u\Delta M(z,u)+M_0(z)$. Then the equation (\ref{eqBMJ}) rewrites to
\begin{equation}\label{eqBMJ-2}
u\Delta M(z,u)+M_0(z) = z Q(z,u, u\Delta M(z,u)+M_0(z),\Delta M(z,u)).
\end{equation}
Hence, the system (\ref{eqBMJsystem}) translates to three equations :
\begin{align*}
u_1(z) f_1(z) + M_0(z) &= zQ(z,u_1(z),u_1(z) f_1(z) + M_0(z), f_1(z)), \\
u_1(z) &=  z u_1(z) Q_{\alpha_0}(z,u_1(z), M_0(z) + u_1(z) f_1(z), f_1(z)) \\
&\quad + zQ_{\alpha_1}(z,u_1(z), u_1(z) f_1(z) + M_0(z), f_1(z))  \\
f_1(z) &= zQ_u(z,u_1(z), u_1(z) f_1(z) + M_0(z), f_1(z) )  \\
&\quad + zf_1(z) Q_{\alpha_0} (z,u_1(z), u_1(z) f_1(z) + M_0(z), f_1(z)).
\end{align*}
If we set $g(z) = u_1(z) f_1(z) + M_0(z)$ then we get a so-called positive polynomial system of equations:
\begin{align*}
g(z) &= zQ(z,u_1(z), g(z), f_1(z)), \\
u_1(z) &=  z u_1(z) Q_{\alpha_0}(z,u_1(z), g(z), f_1(z)) + zQ_{\alpha_1}(z,u_1(z), g(z), f_1(z))  \\
f_1(z) &= zQ_u(z,u_1(z), g(z), f_1(z) )  + zf_1(z) Q_{\alpha_0} (z,u_1(z), g(z), f_1(z)).
\end{align*}
where it holds $f_1(z) = \Delta M(z,u_1(z))$. 
These kinds of systems have been discussed intensively in the literature, see in particular \cite{BD} and
the references there. If this system is non-linear and strongly connected then all solution functions
$g(z), u_1(z), f_1(z)$ have a common dominant square root singularity. Furthermore, if we differentiate 
(\ref{eqBMJ-2}) with respect to $z$ we obtain
\begin{align*}
u(\Delta M)_z(z,u) + M_0'(z) &= Q(z,u,M(z,u),\Delta M(z,u)) 
\\&\;+ zQ_z(z,u,M(z,u),\Delta M(z,u)) 
\\&\;+ z Q_{\alpha_0}(z,u,M(z,u),\Delta M(z,u)) (u (\Delta M)_z(z,u) + M_0'(z))
\\&\; + z Q_{\alpha_1}(z,u,M(z,u),\Delta M(z,u)) (\Delta M)_z(z,u).
\end{align*}
Setting $u=u_1(z)$ and by using the second relation of the above system it follows that all terms
that contain $(\Delta M)_z(z,u)$ cancel and we are able to represent $M_0'(z)$ as
\[
M_0'(z) = \frac{Q(z,u(z), g(z), f_1(z)) + zQ_z(z,u(z), g(z), f_1(z))}{1- zQ_{\alpha_0}(z,u(z), g(z), f_1(z))}.
\]
Since $u(z)(1-zQ_{\alpha_0}) = zQ_{\alpha_1}$ we certainly have $1-zQ_{\alpha_0}\ne 0$ at the singularity.
Hence, $M_0'(z)$ has a square root singularity, too, which implies that $M_0(z)$ has a dominating $3/2$-singularity. (The linear case is even simpler to handle and leads usually to a square root singularity for $M_0(z)$.)

The reason for the relatively simple situation in the case $k=1$ is that it allows to reduce the problem to the solution of a positive polynomial system. Unfortunately, this property is completely lost for $k> 1$.

Let us consider next the case $k=2$. To safe some space we choose the notation $\diff (z,u) = \Delta^{(2)}M(z,u)$ and we set
\[M(z,u) = u^2\diff (z,u)+ uM_1(z) + M_0(z).
\]
The catalytic equation rewrites now to 
\begin{align}\label{eq:seceq}
    u^2\diff(z,u)&+uM_1(z)+M_0(z) \\
    &= zQ\left(z,u, u^2\diff(z,u)+uM_1(z)+M_0(z), u\diff(z,u)+M_1(z), \diff(z,u)\right) \nonumber \\
    &=: R(z,u,\diff(z,u),M_1(z),M_0(z)) \label{eq:R_eq},
\end{align}
where 
\[
R(z,u,y_0,y_1,y_2) = zQ(z,u,u^2y_0+uy_1+y_2,uy_0+y_1,y_0)
\]
is also a polynomial with non-negative coefficients.
%
%
In particular if we translate this equation back to the form (\ref{eqcatBMJ})  we have
\begin{align*}
P(z,u,x_0,x_1,x_2) &= zQ \left(z,u, u^2x_0+ux_1+x_2, ux_0+x_1, x_0\right) - u^2x_0- ux_1- x_2 \\
&= R(z,u,x_0,x_1,x_2) - u^2x_0- ux_1 -x_2.
\end{align*}
Similar to the case $k=1$, the system (\ref{eqBMJsystem}) now rewrites to 
\begin{align}\label{eq:our_system}
    u_i(z)^2f_i(z)+u_i(z)M_1(z)+M_0(z) &= R\left(z,u_i(z), f_i(z),M_1(z),M_0(z)\right),\nonumber\\
    u_i(z)^2 &= R_{y_0}\left(z,u_i(z), f_i(z),M_1(z),M_0(z)\right),
    \\ 
    2u_i(z)f_i(z)+M_1(z) &=  R_u\left(z,u_i(z), f_i(z),M_1(z),M_0(z)\right), \nonumber 
\end{align}
with $i=1$ and $i=2$ for the six indeterminate functions $M_1(z), M_0(z)$, $u_{1,2}(z)$ and $f_{1,2}(z)$ (which correspond to the functions $f_{1,2}(z) = \diff(z,u_{1,2}(z))$). \\

We now state our main results that generalize the results of \cite{DNY2021} to the case $k=2$. However, we make use of the notion of the \emph{dependency graph} of an equation system.

\begin{definition}
Given a positive system of functional equations
\[F_i(z) = Q_i(z,F_1(z),F_2(z),\dots),\qquad  i\geq 1\]
where $Q_i(z,y_0,y_1,\dots)$ are power series with non-negative coefficients in $z$ and $(y_0,y_1,\dots)$, we define the dependency graph as the digraph $G=(V,E)$, where $V= \{F_1,F_2\dots \}$ and the edge $(F_i,F_j) \in E$ if and only if $\frac{\partial}{\partial F_i} Q_j \ne 0$.
We call an equation system strongly connected if its dependency graph is strongly connected.
\end{definition}

Finite strongly connected positive systems have been studied intensely since the 90s (see e.g. \cite{BD}) and the following important property due to the positive coefficients of $Q$ also holds in the infinite case.

\begin{remark}\label{rem:strong}
If the solutions of a positive equation system are algebraic and there is a directed cycle (of finite length) in the dependency graph that connects $F_i(z)$ and $F_j(z)$, then the two functions share a common radius of convergence~$z_0$ and in case $z_0<\infty$, they also show the same singular behaviour at~$z_0$.
\end{remark}

As illustrated by Tutte's planar maps equation system (\ref{eqmapsrec}) a catalytic variable equation is just an infinite system of equations if one considers each equation in the coefficients of $[u^i], i\geq 0$ independently. With this in mind, we may formulate our main result for linear equations as follows.

\begin{theorem}\label{Th1}
Suppose that the polynomial $Q$ in the catalytic equation (\ref{eq:seceq}) is linear in $(\alpha_0, \alpha_1,\alpha_2)$  and has non-negative coefficients. 
Suppose further that the dependency graph of the infinite equation system  
\[ M_i(z) = [u^i] zQ(z,u,M(z,u), \Delta M(z,u), \Delta^{(2)} M(z,u))\]
is strongly connected. 
Then the solutions $M_0(z)$ and $M_1(z)$ have a common radius of convergence $z_0$ and a square root singularity at $z_0$.

Furthermore there exists $d\ge 1$ and a non-empty set $J \subseteq \{0,1,\ldots,d-1\}$ 
of residue classes modulo $d$ and constants $c_j> 0$ such that for $j\in J$
\begin{equation}\label{eqTh1.2-0}
M_n = [z^n]\,M_0(z) = c_j n^{-3/2} z_0^{-n}\left( 1 + O\left( \frac 1n \right) \right), \qquad (n\equiv j \bmod b,\ n\to \infty).
\end{equation}
Furthermore, if $n\equiv j \bmod d$ with $j\not \in J$ then we either have $M_n = 0$ 
or
\[
M_n = O\left( n^{-5/2} z_0^{-n} \right).
\]
A similar statement holds for $[z^n]\,M_1(z)$.
\end{theorem}

\medskip

For the second theorem we will need an extra condition for the term
\begin{align}\label{eq:cond}
T&:= R_{uuu}+(3R_{uuy_0}-6)\frac{2u - R_{uy_0}}{R_{y_0y_0}}+3R_{uy_0y_0}\left(\frac{2u - R_{uy_0}}{R_{y_0y_0}}\right)^2 \\
&+R_{y_0y_0y_0}\left(\frac{2u - R_{uy_0}}{R_{y_0y_0}}\right)^3, \nonumber
\end{align}
that is evaluated at the corresponding values for $z_0$, $u_1(z_0)$, $f_1(z_0)$, $M_0(z_0)$, $M_1(z_0)$,
where $z_0>0$ denotes the (common) dominant singularity.

\begin{theorem}\label{Th2}
Suppose that the polynomial $Q$ in the catalytic equation (\ref{eq:seceq}) is non-linear in $(\alpha_0, \alpha_1,\alpha_2)$ and has non-negative 
coefficients. Suppose further that the dependency graph of the infinite equation system defined by 
\[ M(z,u) = zQ(z,u,M(z,u), \Delta M(z,u), \Delta^{(2)} M(z,u))\]
is strongly connected. 
Then the functions $M_0(z)$ and $M_1(z)$ have a common radius of convergence $z_0>0$.

If $T\ne 0$ then both $M_0(z)$ and $M_1(z)$ have a dominant singularity of type~$3/2$. 
Furthermore, there exists $d\ge 1$ and a non-empty set $J \subseteq \{0,1,\ldots,d-1\}$ 
of residue classes modulo $d$ and constants $c_j> 0$ such that for $j\in J$
\begin{equation}\label{eqTh1.2-0}
M_n = [z^n]\,M_0(z) = c_j n^{-5/2} z_0^{-n}\left( 1 + O\left( \frac 1n \right) \right), \qquad (n\equiv j \bmod b,\ n\to \infty).
\end{equation}
Furthermore, if $n\equiv j \bmod d$ with $j\not \in J$ then we either have $M_n = 0$ 
or
\[
M_n = O\left( n^{-7/2} z_0^{-n} \right).
\]
A similar statement holds for $[z^n]\,M_1(z)$.
\end{theorem}

Indeed, we conjecture that (under the conditions of Theorem~\ref{Th2}) we always have $T \ne 0$. We could not find or computationally construct any example, where $T = 0$. If this conjecture is true, we would have the same kind of universal asymptotic properties as in the (simple) case $k=1$. 
Of course we expect the same phenomenon in general.

\begin{conjecture}
Suppose that $k\ge 1$ and that $M(z,u)$ is the solution of the catalytic equation (\ref{eqBMJgen}),
where $Q$ is a polynomial with non-negative coefficients (that satisfies some proper regularity 
conditions in order to exclude exceptional situations). Then the functions
$M_j(z) = [u^j] M(z,u)$, $0\le j< k$, have a common dominant singularity~$z_0$ that is
of square root type if $Q$ is linear in $(\alpha_0,\ldots,\alpha_{k-1})$ and of type $3/2$ if $Q$ is non-linear in $(\alpha_0,\ldots,\alpha_{k-1})$.
\end{conjecture}

Naturally, one might be interested in technical conditions that guarantee a strongly connected dependency graph. We start with the linear situation, where equation (\ref{eq:seceq}) can be written as
\begin{align}\label{eqMlinear}
M(z,u) &= zQ_0(z,u) + zQ_1(z,u)M(z,u) \\
&\quad + zQ_2(z,u) \Delta M(z,u) + zQ_3(z,u)\Delta^{(2)}M(z,u). \nonumber
\end{align}

In the proof of Theorem \ref{Th1}, we will therefore often simply make use of the technical conditions of Proposition \ref{lem:dep}.
In the following considerations  we will also make use of the so-called curve equation 
$u^2 = R_{y_0}(z,u,\Delta(z,u),M_1(z),M_0(z))$ (see Section~\ref{sec:curve_eq}).
Note that $R_{y_0}$ does not depend on $\Delta(z,u)$, $M_1(z)$,
and $M_0(z)$ in the linear case and is then a polynomial in $z$ and $u$.

\begin{proposition} \label{lem:dep}
    The dependency graph of a linear positive system of the form (\ref{eqMlinear}) is strongly connected if and only if $Q_{1u}(z,u)\neq 0$, the curve equation is not a polynomial in $u^2$ and $u$ does not divide $Q_2(z,u)+Q_3(z,u)$.
\end{proposition}

\begin{proof}
    We will consider the dependency graph of the infinite system 
    \[[u^i] M(z,u) = [u^i] zQ\left(z,u,M(z,u),\diff M(z,u),\diff^{(2)}M(z,u)\right)\]
    where we further write $M(z,u)=\sum_{i\geq 0} M_i(z) u^i$ and first show that, given that $Q_{1u}(z,u)\neq 0$, the curve equation is not a polynomial in $u^2$ and $u^2$ does not divide the curve equation, there is directed cycle on which $M_0(z)$ and $M_1(z)$ lie. By shifting the argument by $i$ to $M_i(z)$ and $M_{i+1}(z)$ the statement then follows immediately.\\
    First, we observe that since $u$ does not divide $Q_2(z,u)+Q_3(z,u)$, there is an edge from $M_i(z)$ to $M_{i-1}(z)$ for all $i\geq 1$ because $Q_2(z,0)\ne 0$ or  there is an edge from $M_i(z)$ to $M_{i-2}(z)$ for $i\geq 2$ in case $Q_3(z,0)\ne 0$.  That is, there exists $j \in \{1,2\}$ such that for all $i\geq j$ there exists an edge in the dependency graph from $M_i(z)$ to $M_{i-j}(z)$.\\
    Further, since $Q_{1u} \neq 0$, there exists a positive integer $k$ such that $[u^k]Q_1(z,u) \ne 0$. One can see from the system that in this case there is an edge from $M_l(z)$ to $M_{l+k}(z)$ for all $l\geq 0$. \\
    Now, we need to distinguish three cases:
    \begin{enumerate}
        \item $j=1$. In this case, we can obviously find the cycle 
        \[M_0(z) \rightarrow M_k(z) \rightarrow M_{k-1} \rightarrow \dots M_1(z) \rightarrow M_0(z)\]
        mutually connecting $M_0(z)$ and $M_1(z)$.
        \item $j=2$ and $k$ is odd. In this case, the cycle is nearly as straight forward as in the first case.
        \[M_0(z) \rightarrow M_k(z) \rightarrow M_{k-2} \rightarrow \dots M_1(z) \rightarrow M_{k+1} \rightarrow \dots \rightarrow M_0(z)\]
        \item $j=2$ and $Q_1(z,u)$ is a polynomial in $u^2$ (otherwise we choose an odd $k$ and apply the second case). Since the curve equation is not a polynomial in $u^2$, there exists $m$ such that $[u^{2m}] Q_2(z,u) \neq 0$ or $[u^{2m+1}] Q_3(z,u) \neq 0$. But then there is an edge from $M_{2m+2}$ to $M_1(z)$. Thus, again we have found a path from $M_0(z)$ to $M_1(z)$ by
        \[M_0(z) \rightarrow M_k(z) \rightarrow \dots \rightarrow M_{sk} \rightarrow \dots \rightarrow M_{2m+2}(z) \rightarrow \dots \rightarrow M_{k-2\left(\frac{k-1}{2}\right)}\]
        where $s$ is the smallest integer such that $sk \geq 2m+2$. Similarly, we can connect $M_1(z)$ to $M_2(z)$ from which there exists an edge to $M_0(z)$ closing the cycle which contains $M_0(z)$ and $M_1(z)$.
    \end{enumerate}
    What is left to show is that the system is not strongly connected in case one of the technical conditions is not satisfied. We treat these cases in Section \ref{sec:deg_linear}. In particular, we see that there is no outgoing edge from $M_0(z)$ in the dependency graph if $Q_{1u}(z,u) = 0$, the dependency graph has two components if the curve equation is an equation in $u^2$ and there is no incoming edge at $M_0(z)$ if $u$ divides $Q_2(z,u)+Q_3(z,u)$.
\end{proof}

For the non-linear case, we make use of the following necessary conditions.
\begin{proposition} \label{prop:nonlintec}
    If the dependency graph of an infinite system of the form
    \[ M(z,u) = zQ\left(z,u,M(z,u), \Delta M(z,u), \Delta^{(2)} M(z,u)\right)\]
    is strongly connected,  then $Q_{\alpha_0}$ depends on $u$ or $\Delta^{(2)} M(z,u)$ and the curve equation is not power series in $u^2$.
\end{proposition}
\begin{proof}
    Suppose $Q_{\alpha_0}$ does not depend on $u$ and $M(z,u)$, then the equation is of the form
    \[ M(z,u) = zQ_0(z)M(z,u) + \tilde Q(z,u,\Delta M(z,u), \Delta^{(2)} M(z,u))\]
    but then there is now outgoing edge from $M_0(z)$ in the dependency graph of the system and the graph is therefore not strongly connected.\\
    Now, assume that the curve equation is a power series in $u^2$. In case $Q$ is linear in $(\alpha_0,\alpha_1,\alpha_2)$ we already have shown that the system is not strongly connected. In case $Q$ is non-linear, the curve equation depends on $\Delta^{(2)}M(z,u)$. Thus, either $[u^{2k}] \Delta^{(2)}M(z,u) =0$ or $[u^{2k+1}] \Delta^{(2)}M(z,u) =0$ for all $k\geq 0$. Therefore, since $M(z,u) \ne 0$, the system is certainly not strongly connected by Remark \ref{rem:strong}.
\end{proof}

\begin{example}\label{ex3}
Let us consider one-dimensional non-negative lattice paths, where we allow steps of the form $\pm 1$ and $\pm 2$.
The generating functions $E_k(z)$ of walks that start at $0$ and end at level $k$ satisfy the system of equations 
\begin{align*}
E_0(z) &= 1 + z(E_1(z) + E_2(z)), \\
E_1(z) &= z(E_0(z) + E_1(z) + E_2(z)), \\
E_k(z) &= z (E_{k-2}(z) + E_{k-1}(z) + E_{k+1}(z) + E_{k+2}(z) ) \qquad (k\ge 2).
\end{align*}
Hence, the generating function $E(z,u) = \sum_{k\ge 0} E_k(z) u^k$ satisfies
\begin{equation}\label{eqEzv}
E(z,u) = 1 + z(u+u^2) E(z,u) + z \frac{E(z,u) - E(z,0)}u + z \frac{E(z,u) - E(z,0) - u E_v(u,0) }{u^2}.
\end{equation}
This is precisely a linear equation of the form (\ref{eqBMJgen}) with $k=2$. 
Theorem~\ref{Th1} applies directly and implies that the generating function $E_0(z) = E(z,0)$ of excursions
has a square root singularity, compare also with \cite{BandFla2002} or with the discussion in Section~\ref{sec6}. 
\end{example}

\begin{example}\label{ex4}
$3$-Constellations are Eulerian maps, where the faces are bi-colored, black faces have valency $3$
whereas white faces have a valency that is a multiple of $3$ (more generally one considers $m$-constellations, see \cite{BM-J}).
The corresponding (catalytic) equation for \mbox{3-constellations} is given by
\begin{align*}
C(z,u) &= 1 + zu C(z,u)^3 + zu(2 C(z,u) + C(z,1)) \frac{C(z,u)- C(z,1)}{u-1}  \\
&+ zu \frac{C(z,u)- C(z,1)- (u-1) C_u(z,1)}{(u-1)^2}
\end{align*}
This catalytic equation is almost of the form, where we can apply Theorem~\ref{Th2} due to the additional 
appearance of $C(z,1)$. However, the polynomial $P$ in (\ref{eq:R_eq}) has still non-negative coefficients.
Thus a slight extension of Theorem~\ref{Th2} applies, {where we require the determinant of $A$ in the calculations of Section \ref{sec3} to be positive at $z_{0}$}. 
Furthermore, (\ref{eq:cond}) is satisfied. 
Consequently, the function $C(z,1)$ has a dominant singularity of type $3/2$, see also the
discussion in Section~\ref{sec6}
\end{example}

Further examples can be found in \cite{BM-J,DNY2021}. It should be also mentioned that Theorems~\ref{Th1} and~\ref{Th2}
can be extended to prove central limit theorems for several parameters that are encoded by an additional variable
(see \cite{DrmPana2013,DrmotaYu2018,DNY2021}).

\begin{theorem}\label{Th3}
Suppose $Q(z,u,w,\alpha_0,\alpha_1,\alpha_2)$ is a polynomial with non-negative coefficients such that
$Q(z,u,1,\alpha_0,\alpha_1,\alpha_2)$ satisfies the assumptions of Theorem~\ref{Th1} or Theorem~\ref{Th2}
and that $T\ne 0$. Furthermore, let $M(z,u,w)$ denote the solution of the catalytic equation
\begin{equation}\label{eqMzuw}
M(z,u,w) = z Q\left(z,u,w, M(z,u,w), \diff (M(z,u,w)) ,\diff^{(2)}(M(z,u,w)) \right)
\end{equation}
where the differences $\diff$ and $\diff^{(2)}$ are taken with respect to $u$, that is, $w$ is considered as
an additional parameter.

Let $j_0\in J$ denote a residue class modulo $d$ (where $d\ge 1$ is givne in 
Theorem~\ref{Th1} or Theorem~\ref{Th2}) for which the asymptotics of the coefficent are of the form
$[z^n]\, M(z,0,1) \sim c_{j_0} n^{-3/2} z_0^{-n}$ or $[z^n]\, M(z,0,1) \sim c_{j_0} n^{-5/2} z_0^{-n}$, respectively,
for some constant $c_{j_0} > 0$. 

If $n\equiv j_0 \bmod d$ let the random variable $X_n$ defined by
\[
\mathbb {P}[X_n = k] = \frac{[z^n w^k]\, M(z,0,w)  }{[z^n]\, M(z,0,1)}.
\]
Then $X_n$ satisfies a central limit theorem of the form
\[
\frac{X_n - \mu n}{\sqrt n} \to N(0,\sigma^2),
\]
where $\mu, \sigma$ are non-negative constants. In particular we have
$\mathbb{E} X_n \sim \mu n$ and $\mathbb{V}{\rm ar} X_n \sim \sigma^2 n$.
\end{theorem}

For example, we can consider (again) one-dimensional non-negative lattice paths, where we
allow steps of the form $\pm 1$ and $\pm 2$, where we additionally count the number of 
$+1$ steps with an additional variable. Then the corresponding equation for
the generating function $E(z,u,w)$ is given by
\begin{align*}
E(z,u,w) &= 1 + z(uw+u^2) E(z,u,w) + z \frac{E(z,u,w) - E(z,0,w)}u \\
&+ z \frac{E(z,u,w) - E(z,0,w) - u E_u(u,0,w) }{u^2}.
\end{align*}
By a direct application of Theorem~\ref{Th3} it follows that the number of $+1$-steps
in paths of length $n$ satisfies a central limit theorem. (It is easy to see, for example,
that $\mu = 1/4$.)

\section{The Curve Equation}\label{sec:curve_eq}

In this section, we will first study general catalytic equations of the type~(\ref{eqBMJgen}) and put special emphasize on what we will call the curve equation. Subsequently, we focus on the case $k=2$ and prove several properties of the solutions to this equation.\\
Our first observation is that $M(z,u)$ is not only uniquely defined as a power series by equation (\ref{eqBMJgen}), but it is also analytic at the origin. We simply consider the equation as a fixed point problem in the sequence space of the coefficients. The factor~$z$ on the right hand side can be chosen small enough such that the map is a contraction and yields a unique solution that is analytic by uniform convergence. Due to the non-negative coefficients of $Q$ (and $R$) it further follows that the solution function has non-negative coefficients.\\

Now, given that we know there is a unique solution 
\[
M(z,u) = u^2\Delta(z,u) + u M_1(z) + M_0(z)
\]
fulfilling the equation, we consider now the equation
\begin{equation}\label{eqcurveequation}
u^k = R_{y_0}\left(z,u, \diff(z,u),M_1(z),M_0(z)\right) =: C(z,u)
\end{equation}
to which we will refer as the \emph{curve equation}. As pointed out in \cite{BM-J}, this equation has exactly $k$ roots $u_1(z),\ldots,u_k(z)$ with $u_i(0) = 0$ and in case they are all distinct, system~(\ref{eqBMJsystem}) is solvable. In this section, we will study the solutions $u_1(z),\dots u_k(z)$ and we will particularly focus on the root with maximal absolute value along the positive real axis which we define to be $u_1(z)$. Note that we exclude the case where $u^k$ divides the curve equation to avoid the trivial case where all $u_i(z) = 0$.

\begin{lemma}\label{le:def_u1}
Let $k\ge 1$ be an integer and suppose that $C(z,u)$ is analytic around $0$ with non-negative coefficients and $z\mid C(z,u)$. Then there exist exactly $k$ Puiseux series $u_1(z),\dots, u_k(z)$ (counted with multiplicities) that satisfy
\[u^k = C(z,u)\]
with $u_i(0)=0$ and if $u^k$ does not divide $C(z,u)$, then there exist a unique real-valued, positive root $u_1(z)$ which takes the maximal absolute value of all $u_i(z)$ for all positive $z$ in a neighbourhood of $0$.\\
Further, there exists $z_0> 0$ such that $u_1(z)$ is analytic for $0< z < z_0$ and not analytic at $z_0$, but 
\[\lim_{z\to z_0-} u_1(z) < \infty.\]
\end{lemma}

\begin{proof}
We may assume that $C(z,0) \ne 0$. Otherwise divide the equation by $u^\ell$, where $\ell$ is the largest integer $< k$ such that $u^\ell \mid C(z,u)$ to avoid trivial solutions $u=0$.\\
As pointed out by Bousquet-M\'elou and Jehanne in \cite{BM-J}, the existence of $k$ roots follows from the Puiseux theorem. \\
Recall that in the Puiseux algorithm, we mark the point $(i,j)$ for each term $a_{ji}u^jz^i$ in the curve equation with $a_{ji} \ne 0$ and the convex hull gives us the Newton polygon. Since $C(z,0) \ne 0$, the polygon intersects the y-axis at some positive integer $\ell$ and the x-axis at $k$ since all terms in $C(z,u)$ appear with a factor $z$. The negative slope of the different segments between these intersection points are the leading exponents of the solutions $u_i(z)$. (These leading exponents are also called valuation of $u_i(z)$ and they are bounded from below by $1/k$.) \\
Of course, plugging a solution $u_i(z) = z^\alpha(c+\tilde{u}_i(z))$ with leading exponent $\alpha$ and a (unique) Puiseux series $\tilde{u}_i(z)$ with positive valuation into the curve equation, the terms with the smallest exponent of $z$ have to cancel in the curve equation. If we choose the slope of a segment which is not incident with $(k,0)$, then these terms only show up in $C(z,z^\alpha(c+\tilde{u}_i(z)))$. Thus, the coefficient $c$ in $z^\alpha(c+\tilde{u}_i(z))$ has to be a zero of a polynomial with positive coefficients and consequently cannot be a positive real value. \\
Therefore we only have to consider solutions with the leading exponent arising from the slope $\ell_0/k_0$ of the segment incident to $(k,0)$. In this case, the coefficient $c$ at the term $z^{-\ell_0/k_0}$ has to satisfy $c^k = P(c)$, where $P$ is a polynomial with non-negative coefficients and degree $<k$ since $z\mid C(z,u)$. It is a direct consequence of the Perron-Frobenius theorem if one considers the companion matrix of the polynomial, that there exists a unique positive root $c_0$ to such an equation. We denote the solution with this coefficient $c_0$ by $u_1(z)$. More precisely,
\[
u_1(z) = \sum_{\ell > \ell'} a_\ell z^{\ell/k'}
\]
with $a_{\ell'} = c_0 > 0$, where $\ell',k'$ are positive integers and $\ell'/k' = -\ell_0/k_0 \geq 1/k$. 
\\
Since $\ell_0/k_0$ is the largest slope of the Newton polygon, the absolute value of $u_1(z)$ will be maximal for small $z>0$. Further, it is not obvious that all coefficients in the expansion are real valued. However, there has to exist a solution that is real valued and positive since for any fixed, positive $u$, it is easy to see that there is exactly one positive value $z>0$ for which the curve equation is satisfied. This continuous curve has to be $u_1(z)$ since it is the only solution with positive leading coefficient.\\

Now, let us show that this solution is analytic in an interval $(0,z_0)$. 
Since $z \mid C(z,u)$ it follows that
\begin{equation}\label{eqderivativeequation}
k u_1(z)^{k-1} > C_u(z,u_1(z))
\end{equation}
if $z> 0$ is sufficiently small. Note also that by the implicit function theorem a condition of 
the form $k u_1(z')^{k-1} \ne C_u(z',u_1(z'))$ for some $z'$ ensures that the function $u_1(z)$
can be analytically continued to a complex neighborhood of $z'$ as long as $C$ is still analytic at $(z',u_1(z'))$. Hence, starting with some sufficiently small $z'> 0$ we can continue $u_1(z)$ analytically along the
positive real axis. \\

The relation (\ref{eqderivativeequation}) also ensures that $u_1(z)$ is strictly increasing
since
\[
u_1'(z) = \frac{C_z(z,u_1(z))}{k u_1(z)^{k-1} - C_u(z,u_1(z))} > 0.
\]
We first observe that $u_1(z)$ cannot stay bounded if $C(z,u)$ is entire and $u_1(z)$ exists for all $z> 0$.
This would contradict the equation (\ref{eqderivativeequation}) since the right hand side would tend to $\infty$ whereas the
left hand side would stay bounded. Thus, if $u_1(z)$ exists for all $z> 0$ then $u_1(z)$ has to tend
to $\infty$ as $z\to \infty$.  However, this is also impossible since by assumption $C(z,u)$ 
has degree $> k$ in $u$. Consequently, there exists $z_0 > 0$ such that $u_1(z)$ exists only 
for $z < z_0$ or for $z \le z_0$. If $C(z,u)$ is not entire, then $z_0$ could also be the radius of convergence of $C(z,u_1(z))$.\\
In both cases, it follows that $u_1(z)$ has to stay bounded
as $z\to z_0$ by the same argument that $C(z,u)$ has degree $>k$ in $u$.
Hence there exists 
\[
\lim_{z\to z_0-} u_1(z)
\]
that is finite and that we will denote by $u_1(z_0)$. 
\end{proof}

In case the initial equation (\ref{eqBMJ-2}) is linear, the curve equation is of course independent of $\Delta^{(2)} M(z,u)$, $M_1(z)$ and $M_0(z)$. Consequently, $C(z,u)$ is polynomial and we can show that $u_1(z)$ has a dominating square root singularity.

\begin{lemma}\label{Le0}
Let $k\ge 1$ be an integer and suppose that $C(z,u)$ is a polynomial with non-negative coefficients,  
with degree $> k$ in $u$,  and with $z\mid C(z,u)$.
Furthermore, let $u_1(z)$ as described in Lemma \ref{le:def_u1}. Then $u_1(z)$ has a square root singularity at $z_0$.
\end{lemma}

For the proof of Lemma~\ref{Le0} we will need the following two general properties (that we state and prove first).
\begin{lemma}\label{Lehelp1}
    Let $D(u) = \sum_{j\geq 0}d_j u^j$ be a power series with non-negative coefficients and at least two positive coefficients. If both $D(u_0) , D'(u_0) < \infty $ for $u_0 > 0$, it holds that
\[
D''(u_0) > \frac{D'(u_0)^2}{D(u_0)} -\frac{D'(u_0)}{u_0}
\]
\end{lemma}

\begin{proof}
We will equivalently show for a fixed positive $u$ that
\[u^2D''(u)D(u) + uD'(u)D(u) > (uD'(u))^2\]
which in turn is satisfied if
\begin{align*}\sum_{k\geq 1} k^2 d_ku^k \cdot \sum_{k\geq 0} d_ku^k &> \left(\sum_{k\geq 1} kd_ku^ k\right)^2\end{align*}
or equivalently
\begin{align*}\sum_{k\geq 1} \left(\sum_{l=0}^k l^2d_ld_{k-l}\right)u^k &> \sum_{k\geq 2} \left(\sum_{l=0}^k l(k-l)d_ld_{k-l}\right)u^k.
\end{align*}
Now, notice that by symmetry, it holds that
\[\sum_{l=0}^k l^2d_ld_{k-l} = \sum_{l=0}^k (k-l)^2d_ld_{k-l} = \sum_{l=0}^k \frac{l^2+(k-l)^2}{2}d_ld_{k-l}\]
and since $a^2+b^2 > 2ab$ for all positive integers with $a \neq b$, the claim follows.
\end{proof}

\begin{lemma}\label{Lehelp2}
Let $F(z,u)$ be a function that is analytic at $(z_0,u_0)$ such that the following conditions hold:
\[
F(z_0,u_0) = 0, \quad F_u(z_0,u_0) = 0, \quad F_z(z_0,u_0) \ne 0, \quad F_{uu}(z_0,u_0) \ne 0.
\]
Then the only local solutions of the equation $F(z,u(z)) = 0$ with $u(z_0) = u_0$ can be represented as
\[
u(z) = g_1(z) \pm g_2(z)\sqrt{z-z_0},
\]
where $g_1(z)$ and $g_2(z)$ are analytic at $z_0$ and satisfy $g_1(z_0) = u_0$ and $g_2(z_0) \ne 0$.
\end{lemma}

\begin{proof}
This is a classical property, see for example \cite[Remark~2.20]{D}.
\end{proof}

Now we can prove Lemma~\ref{Le0}.
\begin{proof}
By Lemma \ref{le:def_u1}, we know that there is a $z_0$ at which $u_1(z)$ becomes singular and takes finite value. Since $C(z,u)$ is a polynomial, we must have
\[
k u_1(z_0)^{k-1} = C_u(z_0,u_1(z_0)).
\]
Otherwise we could find an analytic continuation of $u_1(z)$ around $z_0$ which contradicts
the definition of $z_0$. Another consequence of the definition of $z_0$ is that (\ref{eqderivativeequation}) certainly holds for
all positive $z< z_0$.\\

Now by a direct application of Lemma~\ref{Lehelp1} (where we fix $z_0> 0$ and set $D(u) = C(z_0,u)$) it follows that
\begin{align*}
C_{uu} (z_0,u_1(z_0)) &> \frac{C_u^2}{C} - \frac{C_u}u \\
&= \frac{k^2 u_1(z_0)^{2k-2}} {u_1(z_0)^k} - \frac{k u_1(z_0)^{k-1}}{u_1(z_0)} \\
&= k(k-1)u_1(z_0)^{k-2}.
\end{align*}
Finally we apply Lemma~\ref{Lehelp2} for the function 
\[
F(z,u) = u^k - C(z,u)
\]
at the points $z= z_0$ and $u = u_0 = u_1(z_0)$ and derive
that $u_1(z)$ has a square root singularity at $z_0$.
\end{proof}

From now on, we focus on the case $k=2$ and study both solutions $u_1(z)$ and $u_2(z)$ in further detail.\\
This case is quite advantageous as we can compute both solutions more or less explicitly after the application of the Weierstrass preparation theorem to the curve equation. Clearly, the functions $u_i(z)$ have to be zeros of
\begin{equation}
    \label{eq:weier}
    u^2 - C(z,u) = K(z,u)\left( u^2 + a_1(z) u + a_2(z) \right) = 0
\end{equation}
where $K(z,u)$, $a_1(z), a_2(z)$ are uniquely given analytic functions at $0$ with $a_1(0) = a_2(0) = 0$ and $K(0,0) \neq 0$.  
Hence, we can express
\begin{equation}
u_{1,2}(z) = -\frac{a_1(z)}{2} \pm \sqrt{\frac{a_1(z)^2}{4}-a_2(z)} =: g(z) \pm \sqrt{h(z)}.
\label{eq:g_and_h}
\end{equation}
The functions $g(z)$ and $h(z)$ have some surprising properties which we state in the following Lemmas \ref{Le1} and \ref{Le2}. But first let us give a general idea why these functions are useful.\\

Since $u_1(z)$ is a simple root of the curve equation, $h(z) \ne 0$. So we may split $u$ and all power series in $u$ into two parts: one with a 
factor $\sqrt{h}$ and the other without. For example, we observe that
\[
u_{1,2}^2 = (g \pm \sqrt h)^2 = g^2 + h \pm \sqrt h\, 2g
\]
and in general if $C(z,u)$ is a power series in $z,u$
\begin{align*}
C(z,g \pm \sqrt h) &= \sum_{k\ge 0} C_k(z) (g \pm \sqrt h)^k  \\
&= \sum_{k\ge 0} C_k(z) \sum_{j = 0}^k  {k \choose j} g^{k-j} (\pm 1)^j (\sqrt h)^{j} 
 \\
&= \sum_{k, \ell} C_k(z)  {k \choose 2\ell} g^{k- 2 \ell}  h^{\ell} \pm
\sqrt h  \sum_{k, \ell} C_k(z)  {k \choose 2\ell+1} g^{k- 2 \ell-1}  h^{\ell}.
\end{align*}
In particular we set 
\[
C^+(z,g,h) = \sum_{k, \ell} C_k(z)  {k \choose 2\ell} g^{k- 2 \ell}  h^{\ell}
\]
and 
\[
C^-(z,g,h) = \sum_{k, \ell} C_k(z)  {k \choose 2\ell+1} g^{k- 2\ell-1}  h^{\ell}.
\]
Note that if $C(z,u)$ has just non-negative coefficients the same still holds for $C^+(z,g,h)$ and $C^-(z,g,h)$.
For example we have $(u^2)^+ = g^2 + h$ and $(u^2)^- = 2g$.

In our context $C(z,u) = R_{y_0}(z,u,\diff (z,u),M_1(z),M_0(z))$, where we assume that
$\diff (z,u)$, $M_1(z)$, and $M_0(z))$ are already given. Of course in the linear case
$R_{y_0}$ does not depend on $\diff (z,u)$, $M_1(z)$, and $M_0(z))$ and $R_{y_0}$ is just
a polynomial.

In any case we can rewrite the curve equation to 
 \begin{align*} 
     g^2+h \pm \sqrt h 2g = C^+(z,g,h) \pm \sqrt h\,  C^-(z,g,h) 
 \end{align*} 
 and are lead to consider the system of equations
 \begin{align} \label{eq:gh_sys}
     h = C^+(z,g,h) - g^2, \quad 
     g = \frac{1}{2}C^-(z,g,h).
 \end{align}
It has unique power series solutions which have to be exactly $g$ and $h$ as defined in (\ref{eq:g_and_h}).

The first (and unexpected) property that we prove is that $g(z)$ and $h(z)$ have non-negative coefficients
although the system (\ref{eq:gh_sys}) is not a positive one.

\begin{lemma}\label{Le1}
Suppose that $C(z,u)$ is a power series with non-negative coefficients such that $z$ divides $C(z,u)$ and $u^2$ does not divide $C(z,u)$. 
Furthermore let $u_{1,2}(z) = g(z) \pm \sqrt{h(z)}$ be the two solutions with $u(0) = 0$ of the equation 
$u^2 = C(z,u)$.
Then $g(z)$ and $h(z)$ are power series with $g(0) = h(0) = 0$ and non-negative coefficients. 
\end{lemma}

\begin{proof}
We may assume that $C(z,0)\ne 0$, otherwise $u_2(z) = 0$ and $u_1(z)$ is the unique analytic solution with non-negative coefficients to the equation $u = C(z,u)/u$ (see e.g. \cite{D}). In particular it would hold that $h(z)^2=g(z) = u_1(z)/2$.\\

For the non-degenerate case, we  replace the equation $u_{1,2}(z) ^2 = C(z,u_{1,2}(z))$ by the system
of equations 
\begin{align}
g^2 + h &= C^+(z,g,h),  \label{eqR1} \\
2g &= C^-(z,g,h)   \label{eqR2}
\end{align}
and solve it (uniquely) for power series $g(z), h(z)$ with $g(0) = h(0) = 0$.
(Since $z$ divides $C(z,u)$ the same holds for $ C^{\pm}(z,g,h)$ and consequently
the system (\ref{eqR1})--(\ref{eqR2}) has a unique power series
solution. We just rewrite it to $g = \frac 12 C^-(z,g,h)$, 
$h = C^+(z,g,h)- \frac 14 C^-(z,g,h)^2$ and apply the implicit function theorem or a recursive
argument.)

By differentiating (\ref{eqR1})--(\ref{eqR2}) with respect to $z$ we obtain
\begin{align*}
2gg' + h' &= C_z^+ + C_g^+ g' + C_h^+ h', \\
2g' &= C_z^- + C_g^- g' + C_h^- h'.
\end{align*}
or (if we replace $2g$ by $C^-$)
\begin{align*}
-(C_g^+ - C^-) g' + (1-C_h^+) h' &= C_z^+, \\
(1 - C_g^-/2) g' - (C_h^-/2) h' &= C_z^-/2.
\end{align*}
At this point we mention that 
\[
{k \choose 2\ell} (k-2\ell) - {k \choose 2\ell+1} = {k \choose 2\ell+1} 2\ell \ge 0
\]
which ensures that
\begin{align}
C_g^+(z,g,h) - C^-(z,g,h) &= \sum_{k, \ell} C_k(z)  {k \choose 2\ell} (k-2\ell) g^{k- 2 \ell-1}  h^{\ell} \nonumber\\
&\qquad \qquad - \sum_{k, \ell} C_k(z)  {k \choose 2\ell+1} g^{k- 2\ell-1}  h^{\ell}  \nonumber \\
&= \sum_{k, \ell} C_k(z)  {k \choose 2\ell+1} 2\ell g^{k- 2\ell-1}  h^{\ell} \label{eqdiff}
\end{align}
is a power series with non-negative coefficients.
Hence, we obtain for $g'$ and $h'$
\begin{align}
g' &= \frac {C_z^-(1-C_h^+) + C_z^+ C_h^-} {2 (1-C_h^+)(1-C_g^-/2) - (C_g^+ - C^-)C_h^- } \nonumber \\
&= \frac {C_z^-/2 + \frac{C_z^+ C_h^-}{2(1-C_h^+)}} { 1-C_g^-/2 - \frac{(C_g^+ - C^-)C_h^-/2}{1-C_h^+} } =: T_1(z,g,h), 
\label{eqgder} \\
h' &= \frac { (1-C_g^-/2) C_z^+ + (C_g^+ - C^-) C_z^-/2 } { (1-C_h^+)(1-C_g^-/2) - (C_g^+ - C^-)C_h^-/2 } \nonumber  \\
&= \frac {  C_z^+ +\frac{(C_g^+ - C^-) C_z^-/2}{1-C_g^-/2} } { 1-C_h^+ - \frac{(C_g^+ - C^-)C_h^-/2}{1-C_g^-/2} } =: T_2(z,g,h).
\label{eqder} 
\end{align}
It is clear that $T_1(z,g,h)$ and $T_2(z,g,h)$ are power series in $z,g,h$ with non-negative coefficients. Now we use the relations $g' = T_1(z,g,h)$ and $h' = T_2(z,g,h)$ to prove
inductively that the derivatives $g^{(n)}(0)$ and $h^{(n)}(0)$ are non-negative for all $n\ge 0$.
Clearly for $n= 0$ we have $g(0) = h(0) = 0$. Furthermore for $n=1$ we get 
$g'(0) = T_1(0,0,0) \ge 0$ and $h'(0) = T_2(0,0,0) \ge 0$. 
For $n=2$ we differentiate another time with respect to $z$ and obtain
\begin{align*}
g'' &= T_{1,z}(z,g,h) + T_{1,g}(z,g,h) g' + T_{1,h}(z,g,h) h', \\ 
h'' &= T_{2,z}(z,g,h) + T_{2,g}(z,g,h) g' + T_{2,h}(z,g,h) h'
\end{align*}
and consequently $g''(0)\ge 0$ and $h''(0) \ge 0$. In this way we can proceed further and obtain
$g^{(n)}(0)\ge 0$ and $h^{(n)}(0)\ge 0$ for all $n\ge 0$.\\
\end{proof}

\begin{lemma}\label{Le1b}
Let $C(z,u)$, $g(z),h(z)$ be as in Lemma \ref{Le1}. Suppose $C(z,u)$ is not a power series in $u^2$. Then there exist an integer $d>0$ and integers $d_1,d_2 \ge 0$ such that $g(z) = z^{d_1} \tilde g(z^d)$ and $h(z) = z^{d_2} \tilde h(z^d)$, where $\tilde g(z)$ and $\tilde h(z)$ are power series, where all coefficients are positive.
\end{lemma}
\begin{proof}
We note that $g(0) = h(0) = 0$ and further observe by comparing the corresponding recurrences for the coefficients that the function $\overline g(z)$ and 
$\overline h(z)$ that are given by the system
\begin{align}
\overline g &=T_1(z,\overline g, \overline h), 
\label{eqgder-2} \\
\overline h &=T_2(z,\overline g, \overline h), 
\label{eqhder-2} 
\end{align}
have the same non-zero (and, thus, positive) coefficients as $g(z)$ and $h(z)$, shifted by $z$. Note, that here it is crucial that
$T_1$ and $T_2$ have non-negative coefficients.

By construction the system (\ref{eqgder-2})--(\ref{eqhder-2}) is strongly connected, that is, $T_1$ depends on $\overline h$
and $T_2$ depends on $\overline g$. In particular we can reduce this system to one equation by an elimination process.
That is, we use the first equation to represent $\overline g = G(u,\overline h)$, where $G$ is a power series with
non-negative coefficients that depends on $\overline h$. By substituting the function $G$ in the second equation we
arrive at a single equation for $\overline h$:
\[
\overline h = T_2(z,G(z,\overline h),\overline h) =: H(z,\overline h).
\]
Again the function $H$ has a power series expansion with non-negative coefficients and is non-linear in $\overline h$.
For a single non-linear equation it is known (see \cite{BD}) that the solution function can be written as $\overline h(z) = z^{d_1} \hat h(z^d)$ for some $d\ge 1$, where the power series $\hat h(z)$
has positive coefficients. Similarly we can start with the elimination from the second equation and get a 
single equation for $\overline g$ from which we can deduce that it can be represented as
 $\overline g(z) = z^{d_2} \hat h(z^{d'})$ for some $d'\ge 1$, where the power series $\hat g(z)$
has positive coefficients. However, by inserting these representations in the original system 
(\ref{eqgder-2})--(\ref{eqhder-2}) and by performing substitutions of the form $z \mapsto z e^{2\pi i \ell/d}$ 
and $z \mapsto z e^{2\pi i \ell'/d'}$ is follows that $d=d'$. 

Consequently, $g(z)$ and $h(z)$ can be represented as desired as well.
\end{proof}

Next we determine the singular behavior of $g(z)$ and $h(z)$ if $C(z,u)$ is a polynomial.

\begin{lemma}\label{Le2}
Suppose that $C(z,u)$ is a polynomial with non-negative coefficients,  with degree $\ge 3$ in $u$, $z\mid  C(z,u)$ and $u^2$ does not divide $C(z,u)$.
Denote by $u_{1,2}(z) = g(z) \pm \sqrt{h(z)}$ be the two solutions with $u(0) = 0$ of equation (\ref{eq:weier}).

Then the common convergence radius of $g(z)$ and $h(z)$ is $z_0$ and it is the smallest positive value where $u_1(z)$ is singular. The functions $u_1(z)$, $g(z)$ and $h(z)$ have a common square root singularity at $z_0$, 
whereas $u_2(z)$ is regular for all $0 < z \leq z_0$, that is, the square root singularities of $g(z)$ and $h(z)$ cancel
in the representation $u_2(z) = g(z)-\sqrt{h(z)}$. 
\end{lemma}

\begin{proof} 
First, we point out that by Lemma \ref{Le1}, $g$ and $h$ have non-negative coefficients. Thus, at the (smallest) radius of convergence of these two functions, $u_1(z)$ certainly has to become singular since $h(z) \ne 0$. Conversely, at least one of $g(z)$ and $h(z)$ have to become singular at $z_0$.
\\
For $u_2(z)$, we can state for certain, that it will be regular for all $0<z<z_0$, otherwise $g(z)$ or $h(z)$ would have to become singular as well. Now, let us assume that both $u_1(z)$ and $u_2(z)$ are singular at $z_0$. Then,
\[
u_{i}(z_0)^2 = C(z_0, u_{i}(z_0)), \quad 2 u_{i}(z_0) = C_u(z_0, u_{i}(z_0)),
\]
for both $i=1,2$. This would imply that
\[
2(g_0 + \sqrt{h_0}) = C_u(z_0, g_0 + \sqrt{h_0})\quad \mbox{and} \quad  2(g_0 - \sqrt{h_0}) = C_u(z_0, g_0 - \sqrt{h_0}).
\]
Since 
\[
\frac{C_u(z,g+\sqrt{h}) + C_u(z,g-\sqrt{h})}{2} = C_g^+(z,g,h)
\]
it would follow that
\[
2g_0 = C_g^+(z_0,g_0,h_0)
\]
and therefore
\[
C_g^+(z_0,g_0,h_0) - C^-(z_0,g_0,h_0) = 0.
\]
But as we have seen in the proof of Lemma \ref{Le1}, this is impossible. Hence, $u_2(z)$ is regular at $z_0$ and consequently, 
\[2g(z) = u_1(z)+u_2(z)\] 
and 
\[4h(z) = (u_1(z)-u_2(z))^2\]
share a square root singularity at $z_0$ as $u_1(z_0) > u_2(z_0)$.

\end{proof}

In order to get also asymptotic information on the coefficients we need even more 
precise information on the singularities of $g(z)$ and $h(z)$.

\begin{lemma}\label{Le3}
Suppose that $C(z,u)$ is a polynomial with non-negative coefficients,  with degree $\ge 3$ in $u$,  and with $z\mid  C(z,u)$
and let $g(z) = z^{d_1} \tilde g(z^d)$ and $h(z)=z^{d_2} \tilde h(z^d)$ denote the function that are given by Lemma~\ref{Le2}. Then  $\tilde g(z)$ and $\tilde h(z)$ have a square root singularity at their (common) radius of convergence $z_0^{1/d}$ and there exists $\delta> 0$ such that $\tilde g(z)$ and $\tilde h(z)$ can be analytically continued to the region $\mid z \mid < z_0^{1/d} + \delta$, $\mid z-z_0^{1/d} \mid > 2\delta$.
\end{lemma}
In particular this shows that the only singularities of $g(z)$ and $h(z)$ on the circle of convergence
are $z_0 e^{2\pi i \ell/d}$, $0\le \ell < d$, and that we can apply singularity analysis (see \cite{FO}) to 
get asymptotics for the coefficients of $g(z)$ and $h(z)$. Actually we will not apply this for $g(z)$
and $h(z)$ but for $M_0(z)$, see Lemma~\ref{Le6}.

\begin{proof}
We already know that $g(z)$ and $h(z)$ have a common square root singularity at $z=z_0$.
Consequently the functions $\tilde g(z)$ and $\tilde h(z)$ have a common square root singularity 
at $z = z_0^{1/d}$ which in turn generates the square root singularities of $g(z)$ and $h(z)$ at
$z = z_0 e^{2\pi i \ell/d}$, $1\le \ell < d$.

Of course, the determinant of the Jacobian of the split system (\ref{eqR1})--(\ref{eqR2}) has to become $0$ at the singularity $z_0$, that is
\[
(1-C_h^+)(2-C_g^-) + C_h^-(2g - C_g^+) = (1-C_h^+)(2-C_g^-) + C_h^-(C^- - C_g^+) = 0.
\] 
Equivalently this means that 
\[
D= C_g^-/2 + \frac{(C_g^+ - C^-) C_h^-/2}{1- C_h^+} = 1.
\]
Note that $D = D(z,g,h)$ can be represented as a power series in $z$, $g$, $h$ with
non-negative coefficients. Thus, $z_0$ is the unique positive value for which the equation is satisfied. Now suppose that $\mid z \mid = z_0^{1/d}$ but $z\ne z_0^{1/d}$. Then it follows that
$\mid \tilde g(z) \mid < \tilde g(z_0^{1/d})$ and $\mid \tilde h(z) \mid < \tilde h(z_0^{1/d})$.
Consequently if $\mid z \mid = z_0$ but $z\ne z_0 e^{2\pi i \ell/d}$, $0\le \ell < d$. 
then
\[
\mid g(z) \mid < g(z_0) \quad \mbox{and}\quad \mid h(z) \mid < h(z_0).
\]
Thus, we get
\[
\mid D(z,g(z),h(z)) \mid  \le D(\mid z \mid, \mid g(z_0) \mid, \mid h(z_0) \mid) 
< D(z_0,g(z_0),h(z_0)) = 1
\]
and therefore the Jacobian of the split system is non-zero and so
$g(z)$ and $h(z)$ can be analytically continued. This completes the proof of Lemma~\ref{Le3}.
\end{proof}

\section{Linear Equations}
If  $Q$ is linear in $\alpha_0, \alpha_1$ and $\alpha_2$, we can rewrite (\ref{eq:seceq}) to 
\begin{align*}
M(z,u) &= zQ_0(z,u) + zQ_1(z,u) M(z,u) + zQ_2(z,u) \diff M(z,u) \\
&+ zQ_3(z,u) \diff^{(2)}M(z,u)
\end{align*}
where $Q_0(z,u), Q_1(z,u), Q_2(z,u), Q_3(z,u)$ are polynomials with non-negative coefficients and reformulate it as
\begin{align}
\diff& (z,u)  \left( u^2-z\left(u^2Q_1(z,u)+ u Q_2(z,u)+ Q_3(z,u) \right) \right) \label{eqBMJgen-3-2}\\
&= zQ_0(z,u) + (z Q_1(z,u)-1)M_0(z) + (zuQ_1(z,u)+zQ_2(z,u)-u)M_1(z).  \nonumber
\end{align} 
In this case, the curve equation is a polynomial equation in $u$ and $z$ 
\[u^2-C(z,u) = u^2-z\left(u^2Q_1(z,u)+uQ_2(z,u)+Q_3(z,u)\right) = 0\]
and can be independently solved (and is actually the basic equation of the original kernel method). 
Subsequently, by using the two solutions $u_{1,2}(z)$ of the curve equation we get the following linear system of equations
\begin{align}
M_0(z) 
+   \left(u_1(z) - \frac{zQ_2(z,u_1(z))}{1 -zQ_1(z,u_1(z))}\right) M_1(z) &= \frac{zQ_0(z,u_1(z))}{1 -  zQ_1(z,u_1(z))}, \label{eqM01-1}   \\
M_0(z)
+   \left(u_2(z) - \frac{zQ_2(z,u_2(z))}{1 -zQ_1(z,u_2(z))}\right) M_1(z) &= \frac{zQ_0(z,u_2(z))}{1 -  zQ_1(z,u_2(z))},  \label{eqM01-2}
\end{align}
to calculate $M_0(z)$ and $M_1(z)$. (Of course if these functions are given we can use them
to obtain the full solution function $M(z,u)$.)

\subsection{Proof of Theorem 1}\label{sec2}
The essential step in the proof of Theorem \ref{Th1} is to determine the singular behavior of $M_0(z)$ and $M_1(z)$. But first, we prove a small lemma that justifies the conditions on $Q$. In essence, they ensure that the dependency graph of the infinite system for the functions $M_i(z)$ is strongly connected.

\begin{lemma}\label{Le5}
Suppose that the assumptions of Theorem~\ref{Th1} are satisfied.  
Then $M_0(z)$ and $M_1(z)$ have a square root singularity at $z_0$.

Furthermore there (uniquely) exists an integer $d\ge 1$ 
such that $M_0(z)$ and $M_1(z)$ can be analytically continued to the region
\begin{equation}\label{eqanalyticcontregion}
\mid z \mid < z_0 + \delta, \quad  \mid z - z_0 e^{2\pi i \ell/d} \mid > 2 \delta \quad (0\le \ell < d)
\end{equation}
for any sufficiently small $\delta > 0$. Both functions have a square root singularity at $z_0$.
Moreover, the functions $M_1(z)$ and $M_0(z)$ have square root singularities at the points $z_0 e^{2\pi i \ell/d}$, $1\le \ell < d$.
\end{lemma}

\begin{proof}
We recall that $M_0(z)$ and $M_1(z)$ are given by (\ref{eqM01-1})--(\ref{eqM01-2}) and that 
$u_1(z)$ and $u_2(z)$ are the solutions to the curve equation, where $u_1(z)$ has a square root singularity at
$z_0$, whereas $u_2(z)$ is regular at $z_0$. We recall that (\ref{eqM01-1})--(\ref{eqM01-2}) can be
rewritten as

\begin{equation}\label{eqM0rep}
M_0(z) + \left( u_{1,2}(z) - \frac{zQ_2(z,u_{1,2}(z))}{1- zQ_1(z,u_{1,2}(z)) } \right) M_1(z) 
= \frac{Q_0(z,u_{1,2}(z))}{1- zQ_1(z,u_{1,2}(z))}.
\end{equation}
At this point we rewrite $u_{1,2}(z)$ as $u_{1,2}(z) = g(z) \pm \sqrt{h(z)}$ and split up between the
$+$-part and the $-$-part. In particular we have 
\[
(1)^- = 0 \quad \mbox{and} \quad
\left( u - \frac{zQ_2(z,u)}{1- zQ_1(z,u) } \right)^- = 1 - \left(\frac{zQ_2(z,u)}{1- zQ_1(z,u) } \right)^-,
\]
which leads to 
\[
\left( 1 - \left(\frac{zQ_2(z,u)}{1- zQ_1(z,u) } \right)^- \right) M_1(z) = \left( \frac{Q_0(z,u)}{1- zQ_1(z,u)} \right)^-
\]
Now notice that by our conditions that $Q_{\alpha_0}$ is not a polynomial in $u^2$ the negative part on the right hand side is non-zero.
We then obtain
\[
M_1(z) = \frac{\left( \frac{Q_0(z,u)}{1- zQ_1(z,u)} \right)^-}
{ 1 - \left(\frac{zQ_2(z,u)}{1- zQ_1(z,u) } \right)^-}\Big(z,g(z),h(z)\Big),
\]
where the right hand side depends on $z,g$ and $h$ and has non-negative coefficients. Therefore, it immediately follows that
$M_1(z)$ has a square root singularities at $z_0e^{2\pi i/\ell}$, $0\le \ell < d$, where $d\ge 1$ is given
by Lemma~\ref{Le3}. Furthermore, by applying Lemma~\ref{Le3} (again) it follows that 
$M_1(z)$ can be analytically continued to a region of the form (\ref{eqanalyticcontregion})
for any sufficiently small $\delta > 0$.

Finally, by Lemma \ref{lem:dep}, the same holds for $M_0(z)$.
\end{proof}

\begin{lemma}\label{Le6}
Suppose that the assumptions of Theorem~\ref{Th1} are satisfied and let $d\ge 1$ be
the (unique) integer given in Lemma~\ref{Le5}.  
Then there exists a non-empty set $J \subseteq \{0,1,\ldots,d-1\}$ 
of residue classes modulo $d$ and constants $c_j> 0$ such that for $j\in J$
\begin{equation}\label{eqTh1.2}
M_n = [z^n]\,M_0(z) = c_j n^{-3/2} z_0^{-n}\left( 1 + O\left( \frac 1n \right) \right), \qquad (n\equiv j \bmod b,\ n\to \infty).
\end{equation}
Furthermore, if $n\equiv j \bmod d$ with $j\not \in J$ then we either have $M_n = 0$ 
or
\[
M_n = O\left( n^{-5/2} z_0^{-n} \right).
\]
\end{lemma}

\begin{proof}
Let $J_1$ denote those $\ell \in \{0,\ldots, d-1\}$ for which
$z= z_0 e^{2\pi i \ell/d}$ is also a square root singularity of $M_0(z)$.
Clearly $0\in J_1$ (recall that  $z=z_0$ is a square root singularity of $M_0(z)$).
Suppose that the local expansions around $z= z_0 e^{2\pi i \ell/d}$ (for $\ell \in J_1$)
are given by
\[
M_0(z) = g_\ell(z) - h_\ell(z) \sqrt{ 1- \frac{z}{z_0}  e^{-2\pi i \ell/d} } 
\]
for certain functions $g_\ell,h_\ell$ that are analytic at $z= z_0 e^{2\pi i \ell/d}$.
Then by using standard singularity analysis (see \cite{FO}) together with the properties stated 
in Lemma~\ref{Le5} it follows that the $n$-th coefficient $M_n = [z^n]\,M_0(z)$ is asymptotically
given by 
\[
M_n =  \frac 1{2\sqrt \pi } \sum_{\ell \in J_1} b_\ell  e^{-2\pi i \ell n/d} n^{-3/2} z_0^{-n} 
+ O\left( n^{-5/2} z_0^{-n} \right)
\]
where $b_\ell = h_\ell(z_0 e^{2\pi i \ell/d})$. Hence, if $n \equiv j \bmod d$ then we have
\[
M_n = c_j n^{-3/2} z_0^{-n} + O\left( n^{-5/2} z_0^{-n} \right),
\]
where
\[
c_j = \frac 1{2\sqrt \pi } \sum_{\ell \in J_1} b_\ell  e^{-2\pi i j \ell/d}.
\]
Since the coefficients $M_n$ are non-negative it follows that the numbers $c_j$ are non-negative, too.
Hence Lemma~\ref{Le6} follows by setting $J = \{ j\in \{0,\ldots d-1\} : c_j > 0\}$.
Note that $J$ is non-empty since $z_0$ is a square root singularity of $M_0(z)$.
\end{proof}

Of course this completes the proof of Theorem~\ref{Th1}.

\subsection{Degenerate cases}\label{sec:deg_linear}
Finally we discuss those cases that are not covered by Theorem~\ref{Th1}.
\begin{enumerate}
\item $Q_{1u}(z,u) = 0$.\\
In this case, there is no edge from $M_0(z)$ to any other $M_i(z)$, $i\geq 1$. Thus, we can reduce the recurrence to
\begin{align*}
    \diff M(z,u) = &z\diff Q_0(z,u) + (zQ_1(z,u)+z\diff Q_2(z,u))\diff M(z,u) \\
    &+ (zQ_2(z,u)+z\diff Q_3(z,u))\diff^{2}M(z,u) + zQ_3(z,u)\diff^{(3)}M(z,u)
\end{align*}
to which we may apply Theorem \ref{Th1} or end up again in one of the degenerate cases discussed here. In case $\diff Q_0(z,u) = 0$, $\diff M(z,u) = 0$ obviously solves the equation and $M_0(z) = Q_0(z,u) = Q_0(z)$. In particular, if $Q_2$ is at most linear in $u$ and $Q_3$ at most quadratic in $u$, then the infinite system reduces to a finite linear system with rational solutions.

\item The curve equation $u^2=C(z,u)$ is a polynomial in $u^2$.\\
In this case, we can separate the recurrence into two independent recurrences (with possibly different initial conditions) since in particular $Q_1(z,u), uQ_2(z,u)$ and $Q_3(z,u)$ would be polynomials in $u^2$ and therefore,
\begin{align*}
    M^{\text{even}}(z,u) &= \sum_{k\geq 0 } u^{k}[u^{2k}]\Big(zQ_0(z,u) + zC(z,u)\diff^{(2)} M_{k} + zQ_1(z,u)M_0(z,u)\Big)\\
    &=Q_{0}^{\text{even}}(z,u) + zQ_1(z,\sqrt{u})(u\diff M^{\text{even}}(z,u)+M_0(z)) \\
    &\qquad + z(\sqrt{u}Q_2(z,\sqrt{u})+Q_3(z,\sqrt{u}))\diff M^{\text{even}}(z,u) \\ 
    M^{\text{odd}}(z,u) &= \sum_{k\geq 0} u^{k} [u^{2k+1}]\Big(zQ_0(z,u) + zC(z,u)\diff^{(2)}M(z,u) \\
    &\qquad \qquad \qquad \qquad + z(uQ_1(z,u)+Q_2(z,u))M_1(z,u) \Big)\\
    &= zQ_0^{\text{odd}}(z,u) + + zQ_1(z,\sqrt{u})(u\diff M^{\text{odd}}(z,u)+M_1(z)) \\
    &\qquad + z(\sqrt{u}Q_2(z,\sqrt{u})+Q_3(z,\sqrt{u}))\diff M^{\text{odd}}(z,u) 
\end{align*}
where we define for any series $F^\text{even}(z,u) = \sum_{k\geq 0} u^k [u^{2k}]F(z,u)$ and  $F^\text{odd}(z,u) = \sum_{k\geq 0} u^k [u^{2k+1}]F(z,u)$.\\
Obviously we can apply the theory for first-order discrete differential equations to these two equations, which was treated in \cite{DNY2021}.

\item $u$ divides $Q_2(z,u)+Q_3(z,u)$.\\
In case $u$ divides $Q_2(z,u)+Q_3(z,u)$, $u(z)=0$ is a solution to the curve equation and the original equation~(\ref{eqBMJ-2}) immediately gives a linear equation for $M_0(z)$ if one sets $u=0$. Consequently $M_0(z)$ is a rational function, not depending on any $M_i(z)$, $i\geq 1$.

\item $Q_{1u}(z,u) \ne 0$ and $u$ divides $Q_3(z,u)$ but not $Q_2(z,u)$.\\
In this case we are not really in a degenerate case. Actually all our computations in the proof of Theorem \ref{Th1} still hold. However, there is a more direct way to prove the asymptotic behaviour.\\
Obviously since $u_2(z)=0$ is a solution to the curve equation, we can express $M_0(z)$ in terms of
\begin{align*} M_0(z) &= zQ_0(z,0) + zQ_1(z,0)M_0(z) + zQ_2(z,0) M_1(z)\\
\Longrightarrow M_0(z) &= \frac{zQ_0(z,0)+zQ_1(z,0)M_1(z)}{1-zQ_1(z,0)}
\end{align*}
The second solution $u_1(z)$ to the curve equation satisfies
\[u = zuQ_1(z,u)+zQ_2(z,u) +z\frac{Q_3(z,u)}{u}\]
to which we can apply standard theory. By the assumption $Q_{1u}(z,u) \ne 0$, the solution $u_1(z)$ has a square root singularity at some positive $z_0$ (see \cite{Bender}) and we can compute $M_1(z)$ by plugging $u_1(z)$ into the original equation, subtracting the equation for $u_2(z)=0$ and dividing by $u_1(z)$ to obtain
\begin{align*} M_1(z) = &\;z\diff Q_0(z,u_1(z)) +zQ_1(z,u_1(z))M_1(z) \\&\;+ z\diff Q_1(z,u_1(z)) M_0(z) + z\diff Q_2(z,u_1(z)) M_1(z)\end{align*}
From here it is easy to see, that $M_1(z)$ has a square root singularity, because
\begin{align*}
    M_1(z) = &\; \frac{z\diff Q_0(z,u(z)) + z\diff Q_1(z,u(z))\frac{zQ_0(z,0)}{1-zQ_1(z,0)}}{1-zQ_1(z,u(z)) - z\diff Q_1(z,u(z)) \frac{zQ_1(z,0)}{1-zQ_1(z,0)} - z\diff Q_2(z,u(z))}.
\end{align*}
\end{enumerate}
Of course this square root singularity translates to $M_0(z)$ and the description of the coefficients follow from the transfer lemma.

\section{Non-Linear Equations}

The main difference to the linear case is that the curve equation is not independent of $\diff(z,u)$. Thus, the singularity of $u_1(z)$ at $z_0$ may be induced by the convergence radius of $\diff(z,u)$, $M_1(z)$ and/or $M_0(z)$. Therefore we consider the full system~(\ref{eq:our_system}) in the order
\begin{align}
u_1^2 f_1 + u_1M_1 + M_0 - R(z,u_1,f_1,M_1,M_0) & = 0, \nonumber \\
u_2^2 f_2 + u_2M_1 + M_0 - R(z,u_2,f_2,M_1,M_0) & = 0 \nonumber  \\
u_1^2 - R_{y_0}(z,u_1,f_1,M_1,M_0) &= 0,   \label{eq:our_system-2} \\
2u_1f_1 + M_1 - R_u(z,u_1,f_1,M_1,M_0) &=0,\nonumber  \\
u_2^2 - R_{y_0}(z,u_2,f_2,M_1,M_0) &= 0,  \nonumber \\
2 u_2 f_2 + M_1 - R_u(z,u_2,f_2,M_1,M_0) &=0. \nonumber 
\end{align}
and its Jacobian matrix at $z_0$ with the goal to determine the types of singularities of the involved functions. 

\subsection{Proof of Theorem~\ref{Th2}}\label{sec3}

Analogously to the linear case, the proof of Theorem~\ref{Th2} will mostly concern the singularities of $u_1(z), g(z)$ and $h(z)$. By Lemma~\ref{Le1}, we already know that $u_1(z_0) < \infty$. Due to the positive coefficients of $R_{y_0}$, $z_0$ is bounded by the radius of convergence of $M_0(z), M_1(z)$ and $f_1(z)$ and further they also take finite value at $z_0$. Another reason to focus on $u_1(z)$ is the fact that $\mid u_2(z)\mid \le u_1(z)$, that $u_2(z)$ is regular for all $0<z<z_0$,  and as $\diff(z,u)$ has non-negative coefficients, $f_2(z)=\diff(z,u_2(z))$ is as well regular in this range.\\

First, we analyze the Jacobian matrix of system~(\ref{eq:our_system-2}) at $z_0$ and show that the principal submatrix corresponding to only five of the six equations is invertible -- we will leave out the third equation. That means we can compute functions $\bar{f}_1(z,u_1), m_1(z,u_1), m_0(z,u_1), \bar{u}_2(z,u_1), \bar{f}_2(z,u_1)$ that are analytic at $z_0, u_1(z_0)$ by the implicit function theorem. Of course, evaluated at $u_1(z)$ these functions are exactly the unknown functions locally around $z_0$. Given these functions, we further compute their partial derivatives with respect to $u_1$ such that we can verify similar conditions as in Lemma~\ref{Lehelp2} in case $T\ne 0$.\\
In a last step, we consider the split system and study the singular behavior of $g(z)$, $h(z)$ and subsequently $M_0(z)$ and $M_1(z)$ on the circle of convergence.\\

In order to shorten notation, we set \vspace{-3mm}
\[ u = u_1(z),\qquad \overline u = u_2(z)\vspace{-3mm}\]
    \[f = f_1(z) = \diff(z,u_1(z)),\qquad  \overline f = f_2(z) = \diff(z,u_2(z)).\vspace{1mm}\]

Similarly, if an expression like \mbox{$R(z,u,\diff(z,u),M_1(z),$ $M_0(z))$} 
is evaluated along $u_1(z)$ (and $f_1(z)=\diff(z,u_1(z))$), we just write $R$. If the expression is evaluated along $u_2(z)$, we will write $\overline{R}$.\\
So the Jacobian matrix of system~(\ref{eq:our_system-2}) with respect to $(M_0,M_1, u, f, \overline u, \overline f)$ is given by
\[\left(\begin{matrix}
    A &0 &0\\
    C_1 & B_1 & 0\\
    C_2 &0 & B_2
\end{matrix}\right)\]
where $A,C_1,C_2,B_1,B_2$ are $2\times 2$ matrices and equal
\begin{align*}
   A = \left(
\begin{matrix}
1-R_{y_2} & u-R_{y_1} \vspace{2mm} \\
\overline{1 - R_{y_2}} & \overline{u-R_{y_1}}
\end{matrix}
\right), \quad 
B_1 &= \left(
\begin{matrix}
2u - R_{uy_0} & -R_{y_0y_0} \vspace{2mm} \\ 
2f - R_{u u} & 2u - R_{uy_0}
\end{matrix}
\right), \\ 
B_2 &= \left(
\begin{matrix}
\overline{2u - R_{uy_0}} & \overline{-R_{y_0y_0}} \vspace{2mm} \\
\overline{2 f - R_{u u}} & \overline{2u - R_{uy_0}}
\end{matrix}
\right).
\end{align*}
The zero entries in the first two rows of the Jacobian correspond to the partial derivatives of the original equation with respect to $u$ and $y_0$ which are satisfied as they already appear in~(\ref{eq:our_system-2}). The other zeroes appear because, for example, $\overline u$ does not appear in the equation $u^2 - R_{y_0} = 0$ and vice versa.

Clearly the determinant of the Jacobian factors into the three determinants 
of the submatrices and not very surprisingly, the determinants of $B_1$ and $B_2$ (nearly) equal the partial derivative with respect to $u$ of the curve equation.

\begin{lemma}\label{lem:B_1B_2help}
    Let $z_0$ be the smallest positive $z$, where $u_1(z)$ is singular. Then for all $0<z< z_0$ it holds that $\det B_1, \det B_2 > 0$ and
        \[2u_1(z) - C_u(z,u_1(z)) = \sqrt{\det B_1},\quad  2u_2(z) - C_u(z,u_2(z))=-\sqrt{\det B_2}.\]
\end{lemma}

\begin{proof} If we differentiate the original equation~(\ref{eqBMJ-2}) twice with respect to $u$, we obtain
\begin{equation}\label{eq:sec_der_u}
    \left(u^2-R_{y_0}\right) \partial_u^2 \diff + 2\left(2u-R_{uy_0}\right)\partial_u\diff +2\diff = R_{uu} + R_{y_0y_0}\left(\partial_u \diff\right)^2
\end{equation}
Thus, along $u_1(z)$ and $u_2(z)$ it holds that
\begin{equation}
        \label{eq:der_u}
        \diff_u = \frac{2u - R_{uy_0} \pm \sqrt{\det B_1}}{R_{y_0y_0}},\qquad \overline{\diff}_u = \overline{\frac{2u - R_{uy_0} \pm \sqrt{\det B_2}}{R_{y_0y_0}}}.
\end{equation}
If we plug these expressions into the curve equation, we end up with
\[2u - C_u(z,u) = 2u - R_{uy_0} - R_{y_0y_0}\frac{2u - R_{uy_0} \pm \sqrt{\det B_1}}{R_{y_0y_0} } = \mp \sqrt{\det B_1}\]
and the same for $u_2(z)$.
Further note that for $z$ small, plugging the Puiseux expansion of $u_{1,2}(z)$ into the derivative of the curve equation shows that it will take positive value along $u_1(z)$ and negative along $u_2(z)$ and therefore the correct choices of signs are $-\sqrt{\det B_1}$ and $+\sqrt{\det B_2}$ in the expressions of $\diff_u$ and $\overline{\diff_u}$. Since both $u_1(z)$ and $u_2(z)$ are real valued, we also know that $\det B_1,\det B_2$ are positive for $0<z<z_0$.
\end{proof}

Next we show that $\det A$ is never $0$ in $(0,z_0)$ and that the submatrix $B_2$ which corresponds to the equations for $u_2(z)$ and $f_2(z)$ is invertible if~\mbox{$\det B_1=0$}.

\begin{lemma} 
Let $z_0$ be the smallest positive $z$, where $u_1(z)$ is singular. Then, the determinants $\det A, \det B_1, \det B_2$ evaluated at $z_0$ satisfy 
$\det A \neq 0$, $\det B_1 =0$, and $\det B_2 \neq 0$.
\end{lemma}

\begin{proof}
We first note that $R_{y_0} = zu^2Q_{\alpha_0} + zuQ_{\alpha_1} + zQ_{\alpha_2}$, 
$R_{y_1} = u zQ_{\alpha_0} + zQ_{\alpha_1}$, and $R_{y_2} = zQ_{\alpha_0}$. Clearly, by comparing the expressions to the curve equation $u^2 = R_{y_0}$, it holds that
\[
u^2 (1- R_{y_2}) = zuQ_{\alpha_1} + zQ_{\alpha_2} \quad \mbox{and}\quad
u (u- R_{y_1}) =  zQ_{\alpha_2}.
\]
Hence, we know along $u=u_1(z)$ 
\[
1- R_{y_2} > 0 \quad \mbox{and} \quad  u- R_{y_1} > 0.
\]
Furthermore, since $u_1(z) > 0 > u_2(z)$ and by Proposition \ref{prop:nonlintec}, it holds that $1- R_{y_2} < \overline {1- R_{y_2}}$ since any expression with non-negative coefficients depending on $u$ (or $\diff(z,u)$) evaluated at $u_2(z)$ will take smaller absolute value than evaluated at $u_1(z)$.
If $\det A =0$, it would also have to hold that
    $u-R_{y_1} < \overline{u-R_{y_1}}$. However, the derivative of equation (\ref{eq:seceq}) with respect to $z$ 
\[
u^2 \diff_z + u M_1' + M_0' = R_z + R_{y_0} \diff_z + R_{y_1} M_1' + R_{y_2} M_0'.
\]
can be rewritten in terms of
\begin{equation} \label{eq:matrixA}\left(
\begin{matrix}
1-R_{y_2} & u-R_{y_1} \vspace{2mm} \\
\overline{1 - R_{y_2}} & \overline{u-R_{y_1}}
\end{matrix}\right) \left(\begin{matrix}
M_0' \vspace{2mm} \\ M_1'
\end{matrix}\right) = \left(\begin{matrix}
R_z \vspace{2mm}  \\ \overline{R_z}
\end{matrix}\right)\end{equation}
as the terms with factor $\diff_z$ add up to the curve equation again and consequently cancel along $u_1(z)$ and $u_2(z)$.\\
But $M_1(z)$ and $M_2(z)$ have non-negative coefficients and the equation (\ref{eq:matrixA}) would imply that, even if $M_1'(z),M_0'(z)$ diverge at $z_0$, $R_z < \overline{R_z}$ in a neighbourhood of $z_0$. This is certainly wrong because any expression with positive coefficients evaluated at $u_1(z)$ is at least as large as the expression at $u_2(z)$. Thus it follows $\det A \ne 0$ and by inverting $A$ and continuity of all functions it also holds that $M_1'(z_0),M_0'(z_0) < \infty$.\\

Next in order to prove that $\det B_1$ and $\det B_2$ cannot be both $0$ we consider the curve equation 
$u^2 = C(z,u) = R_{y_0}(z,u,\diff(z,u),M_1(z),M_0(z))$ and observe that $C_u(z,u)$ is given by
\begin{align*}
    C_u(z,u) = R_{uy_0}&(z,u,\diff(z,u),M_1(z),M_0(z)) \\&+ R_{y_0,y_0}(z,u,\diff(z,u),M_1(z),M_0(z)) \diff_u(z,u).
\end{align*}

Now let us assume that $\det B_1=\det B_2=0$. By Lemma~\ref{lem:B_1B_2help}, this would mean that 
\[2u_i(z_0) = C_u(z_0,u_i(z_0)), \quad i=1,2.\]
However, remember that for $0<z<z_0$ it holds that $2g = C_u^+$. By continuity, this extends to $z=z_0$ and we have already seen in the proof of Lemma~\ref{Le2} that $C_u^+-C^-$ has non-negative coefficients which contradicts the assumption that
\[2u_1(z_0)+2u_2(z_0) - C_u(z_0,u_1(z_0)) - C_u(z_0,u_2(z_0)) = 0.\]

This is also the reason why it must hold $\det B_1=0, \det B_2\ne 0$ and not the other way round since for $0<z<z_0$,
\[\sqrt{\det B_1} + (-\sqrt{\det B_2}) = 2(2g-C^{+}_u) = 2(C^--C^+_u) <0.\]
\end{proof}

Intuitively, this lemma tells us that 'the origin of the singularities' of the involved functions are $u_1(z)$ and $f_1(z)$.

\begin{lemma}
    Let $u_1(z),u_2(z),f_1(z),f_2(z),M_0(z), M_1(z)$ be the solutions to system (\ref{eq:our_system-2}) and let $z_0$ be the smallest positive $z$ where $u_1(z)$ is singular. Then $u_1'(z)$ and $f_1'(z)$ diverge at $z_0$ whereas the derivatives of the other functions take finite value at $z_0$.
\end{lemma}

\begin{proof}
Using the implicit function theorem, we can compute the partial derivatives locally around a point $z'<z_0$ as follows. The function
\begin{align*}
F: (z,u_1,u_2,f_1,f_2M_1,M_0) \mapsto \left(\begin{matrix}u_1^2 f_1 + u_1M_1 + M_0 - R(z,u_1,f_1,M_1,M_0)\\
u_2^2 f_2 + u_2M_1 + M_0 - R(z,u_2,f_2,M_1,M_0)\\
2u_1f_1 + M_1 - R_u(z,u_1,f_1,M_1,M_0) \\
u_2^2 - R_{y_0}(z,u_2,f_2,M_1,M_0)\\
2 u_2 f_2 + M_1 - R_u(z,u_2,f_2,M_1,M_0) 
\end{matrix}\right)
\end{align*}
is analytic at $(z', u_1(z'), u_2(z'), f_1(z'), f_2(z'), M_1(z'),M_0(z'))$ and satisfies 
\[F\big(z, u_1(z), u_2(z), f_1(z), f_2(z), M_1(z),M_0(z)\big) = 0\] 
for analytic functions $u_1(z), u_2(z), f_1(z), f_2(z), M_1(z), M_0(z)$ locally around $z'$.
Now, if we compute the derivatives with respect to $z$ of the left hand side of the equation, they are of course all~$0$ and we obtain for example
\begin{align*}
    \left(\begin{matrix}
    A &0 &0 \\
    C_1 & B_1 &0\\
    C_2&0 & B_2
    \end{matrix}\right) \partial_z\left(\begin{matrix}
    M_0 \\ M_1 \\ u_1 \\ f_1 \\ u_2 \\ f_2
    \end{matrix}\right) (z' ) = \left(\begin{matrix}
    -R_z \\ -\overline{R_z} \\ -R_{zy_0} \\ -R_{zu} \\ -\overline{R_{zy_0}} \\ -\overline{R_{zu}}
    \end{matrix}\right)(z')
\end{align*}
Since the Jacobian matrix of the function is invertible, it therefore holds
\begin{align*}
    \partial_z\left(\begin{matrix}
    M_0 \\ M_1 \\ u_1 \\ f_1 \\ u_2 \\ f_2
    \end{matrix}\right) (z' ) = \left(\begin{matrix}
    A^{-1} &0 &0 \\
    -B_1^{-1}C_1A^{-1} & B_1^{-1} &0\\
    -B_2^{-1}C_2A^{-1} &0 & B_2^{-1}
    \end{matrix}\right)\left(\begin{matrix}
    -R_z \\ -\overline{R_z} \\ -R_{zy_0} \\ -R_{zu} \\ -\overline{R_{zy_0}} \\ -\overline{R_{zu}}
    \end{matrix}\right)(z')
\end{align*}
for $z'<z_0$. However, by continuity $M_0'(z), M_1'(z), u_2'(z), f_2'(z)$ will take finite value at $z_0$ because $A$ and $B_2$ are invertible at $z_0$. Hence the types of the singularities of $M_0(z), M_1(z), u_2(z)$ and $f_2(z)$ (if the functions are singular at $z_0$) are larger than $1$. The values of $u_1'(z), f_1'(z)$ on the other hand are the limit of
\begin{align*}
    &\left(\begin{matrix}
    -B_1^{-1}C_1A^{-1}, & B_1^{-1}
    \end{matrix}\right)\left(\begin{matrix}
    -R_z \\ -\overline{R_z} \\ -R_{zy_0} \\ -R_{zu} 
    \end{matrix}\right) \\&= \frac{1}{\det B_1}\left(\mbox{adj}(B_1)\left(\begin{matrix}
    -R_{zy_0} \\ -R_{zu}
    \end{matrix}\right) -\frac{1}{\det A}\mbox{adj}(B_1)C_1 \mbox{adj}(A)\left(\begin{matrix}
    -R_z \\ -\overline{R_z}
    \end{matrix}\right) \right) 
\end{align*}
where adj(A) and adj($B_1$) denote the adjoints of the matrices $A$ and $B_1$ respectively. By (\ref{eq:matrixA}), the above equals
\begin{align*}
    \left(\begin{matrix}
    -B_1^{-1}C_1A^{-1}, & B_1^{-1}
    \end{matrix}\right)&\left(\begin{matrix}
    -R_z \\ -\overline{R_z} \\ -R_{zy_0} \\ -R_{zu} 
    \end{matrix}\right)
    = \frac{1}{\det B_1}\mbox{adj}(B_1)\left(\left(\begin{matrix}
    -R_{zy_0} \\ -R_{zu}
    \end{matrix}\right) - C_1 \left(\begin{matrix}
    -M_0'(z) \\ -M_1'(z)
    \end{matrix}\right) \right) \\
    & = \frac{1}{\det B_1}\mbox{adj}(B_1) \left(\begin{matrix}
    -R_{zy_0}- R_{y_0y_2}M'_0(z) - R_{y_0y_1}M_1'(z) \\ -R_{zu} - R_{uy_2}M'_0(z) - (1-R_{uy_1})M'_1(z)
    \end{matrix}\right)
\end{align*}
If we divide by $R_{y_0y_0}$ the first entry is equal to
\begin{align*} \frac{1}{\det B_1} \Bigg(\frac{2u-R_{uy_0}}{R_{y_0y_0}}&\left(-R_{zy_0}- R_{y_0y_2}M'_0(z) - R_{y_0y_1}M_1'(z)\right) \\
&-R_{zu} - R_{uy_2}M'_0(z) - (1-R_{uy_1})M'_1(z)\Bigg)\end{align*}
Now let us do a similar trick as in the computation of $\diff_u(z,u)$. We consider 
\begin{align} \label{eq:partialuz}
    \partial_z \partial_u &\left( u^2\diff(z,u) + uM_1(z)  +M_0(z) - R(z,u,\diff(z,u),M_1(z),M_0(z))\right)\\
    &= \left(u^2 - R_{y_0}\right) \diff_{zu}(z,u) + \left(2u-R_{uy_0}-R_{y_0y_0}\diff_u(z,u) \right) \diff_z(z,u) \nonumber \\
    &\quad + M_1'(z) - R_{zu} - R_{zy_0}\diff_u(z,u)
     - R_{uy_1}M_1'(z) - R_{uy_2}M_0'(z)  \nonumber\\
    &\quad - R_{y_0y_1}\diff_u(z,u)M_1'(z) - R_{y_0y_2} \diff_u(z,u)M_0'(z)  \nonumber \\
    &= (1-R_{uy_1} - R_{y_0y_1}\diff_u(z,u))M_1'(z) - (R_{uy_2}+R_{y_0y_2}\diff_u(z,u))M_0'(z) \nonumber  \\
    &\quad- R_{zu} - R_{zy_0}\diff_u(z,u) +\left(2u-R_{uy_0}-R_{y_0y_0}\diff_u(z,u) \right) \diff_z(z,u)= 0 \nonumber
\end{align}
The terms with factor $\diff_{zu}(z,u)$ add up to the curve equation and cancel, while the ones with factor $\diff_z(z,u)$ only add up to $0$ at $(z_0,u_0)$ since equation (\ref{eq:second}) holds. At $z_0$, where $(2u-R_{uy_0})/R_{y_0y_0} = \diff_u(z,u)$, we therefore obtain
\begin{align*}\lim_{z\rightarrow z_0} \frac{u_1'(z)}{R_{y_0y_0}} &= \lim_{z\rightarrow z_0} \frac{\left(2u-R_{uy_0}-R_{y_0y_0}\diff_u(z,u) \right) \diff_z(z,u)}{\det B_1} \\
&= \lim_{z\rightarrow z_0} \frac{\sqrt{\det B_1} \diff_z(z,u)}{\det B_1} \end{align*}
As $\diff_z(z,u)>0$ for all positive $z,u$, $u_1'(z)$ and consequently also $f_1'(z)$ diverge at $z_0$ and their type of singularity has to be in the interval $(0,1)$.
\end{proof}

Equation (\ref{eq:der_u}) concerning the partial derivative of $\diff(z,u)$ further tells us that $2u > R_{uy_0}$ at $z_0$ since $\diff(z,u)$ has non-negative coefficients and $u_1(z)>0$ for $z>0$. This means that the following submatrix of the Jacobian of system~(\ref{eq:our_system}) is invertible and its inverse equals
\[\left(\begin{matrix}
    A &0 &0\\
    C_1 & 2u - R_{uy_0} & 0\\
    C_2 &0 & B_2
\end{matrix}\right)^{-1} = \left(\begin{matrix}
    A^{-1} &0 &0\\
    D & (2u - R_{uy_0})^{-1} & 0\\
    -B_2^{-1}C_2A^{-1} &0 & B_2^{-1}
\end{matrix}\right)\]
where $C_1,C_2$ are generally non-zero matrices that contain the partial derivatives with respect to $M_0$ and $M_1$ of the third, fifth and sixth equation of system (\ref{eq:our_system-2}) respectively and $D$ is a $1\times 2$ matrix.
The implicit function theorem therefore yields analytic functions $m_0(z,u_1)$, $m_1(z,u_1)$, $\tilde u_2(z,u_1)$, 
$\tilde f_1(z,u_1)$, and $\tilde f_2(z,u_1)$ with $m_0(z_0,u_1(z_0)) = M_0(z_0)$, 
$m_1(z_0,u_1(z_0)) = M_1(z_0)$, $\tilde u_2(z_0,u_1(z_0)) = u_2(z_0)$, 
$\tilde f_1(z_0,u_1(z_0)) = f_1(z_0)$, and $\tilde f_2(z_0,u_1(z_0)) = f_2(z_0)$ and we can compute their partial derivatives with respect to $u_1$ (evaluated at $u_0 = u_1(z_0)$).

\begin{lemma}\label{lem:derivatives}
 Let $m_0(z,u_1)$, $m_1(z,u_1)$ and $\tilde f_1(z,u_1)$ be as above. Then for $(z_0,u_0)$ with $u_0=u_1(z_0)$ it holds that
\begin{align*}
 \partial_{u_1} m_0(z_0,u_0) &= \partial_{u_1} m_1(z_0,u_0)= \partial_{u_1} \tilde f_2(z_0,u_0)= \partial_{u_1} \tilde u_2(z_0,u_0) = 0\\
 \partial_{u_1} \tilde f_1(z_0,u_0)  &= \partial_{u}\diff(z_0,u_0)\\
 \partial^2_{u_1}m_{0}(z_0,u_0) &= \partial^2_{u_1}m_{1}(z_0,u_0) = 0\\ 
 \partial^2_{u_1}\tilde f_{1}(z_0,u_0) &= \left(\frac{R_{uuu}
+2(R_{uuy_0}-2)\partial_{u_1}\tilde f_{1}+R_{uy_0y_0}(\partial_{u_1}\tilde f_{1})^2}{2u_1-R_{uy_0}}\right)(z_0,u_0)
\end{align*}
\end{lemma}
\begin{proof}
As we discussed earlier, we compute the partial derivatives with respect to $u_1$ using the implicit function theorem. We obtain for example
\begin{equation*} \label{eqderivativeu1}
   \left(\begin{matrix}
    A &0 &0\\
    C_1 & B_1 & 0\\
    C_2 &0 & B_2
\end{matrix}\right) \partial_{u_1} \left(\begin{matrix}
        m_0(z,u_1)\\
        m_1(z,u_1)\\
        \tilde f_1(z,u_1)\\
        \tilde u_2(z,u_1)\\
        \tilde f_2(z,u_1)
    \end{matrix} \right) = - \left(\begin{matrix}
    0\\
    0\\
    2\tilde f_1 - R_{uu}\\
    0\\
    0
\end{matrix} \right)
\end{equation*}
or equivalently,
\begin{equation} \label{eqderivativeu1}
   \partial_{u_1} \left(\begin{matrix}
        m_0(z,u_1)\\
        m_1(z,u_1)\\
        \tilde f_1(z,u_1)\\
        \tilde u_2(z,u_1)\\
        \tilde f_2(z,u_1)
    \end{matrix} \right) = - \left(\begin{matrix}
    A^{-1} &0 &0\\
    D & (2u_1 - R_{uy_0})^{-1} & 0\\
    -B_2^{-1}C_2A^{-1} &0 & B_2^{-1}
\end{matrix}\right) \left(\begin{matrix}
    0\\
    0\\
    2\tilde f_1 - R_{uu}\\
    0\\
    0
\end{matrix} \right)
\end{equation}
Clearly, if we evaluate these expressions at $(z_0,u_0)$, we compute the first two identities of the lemma. That is
\begin{equation*}
 \partial_{u_1} m_0(z_0,u_0) = \partial_{u_1} m_1(z_0,u_0) = \partial_{u_1} \tilde f_2(z_0,u_0)= \partial_{u_1} \tilde u_2(z_0,u_0)= 0
\end{equation*}
and since $\det B_1 = 0$ it further holds that
\begin{equation}\label{eq:partialf_1} \partial_{u_1} \tilde f_1(z_0,u_0) = \left(-\frac{2\tilde f_1- R_{uu}}{2u_1-R_{uy_0}}\right)(z_0,u_0) = \left(\frac{2u_1-R_{uy_0}}{R_{y_{0}y_{0}}}\right)(z_0,u_0) = \diff_u(z_0,u_0).\end{equation}

The same strategy works of course for higher derivatives. If we only consider the submatrix for the relevant derivatives, we obtain that
    \begin{align*}
   \partial^2_{u_1} &\left(\begin{matrix}
        m_0\\
        m_1\\
        \tilde f_1
    \end{matrix} \right) \\&= \left(\begin{matrix}
    A^{-1} &0 \\
    D & (2u_1 - R_{uy_0})^{-1}
\end{matrix}\right) 
   \cdot \left(\begin{matrix}
    R_{uu}+(2R_{uy_0}-4u_1)(\partial_{u_1} \tilde f_{1}) + R_{y_0y_0} (\partial_{u_1}\tilde f_{1})^2-2\tilde f_1 \\
    0\\
    R_{uuu}+2R_{uuy_0}\partial_{u_1}\tilde f_{1}+R_{uy_0y_0}(\partial_{u_1}\tilde f_{1})^2-4\partial_{u_1}\tilde f_{1}
\end{matrix} \right)
\end{align*}
Now by dividing the first entry in the right vector above by $R_{y_0y_0}$, it is easy to see that at $(z_0,u_0)$ it is equal to $0$ because by equation (\ref{eq:partialf_1}),
\[
\left(\frac{R_{uu}-2\tilde f_1}{2u_1-R_{uy_0}}\frac{2u_1-R_{uy_0}}{R_{y_0y_0}}-2\frac{2u_1-R_{uy_0}}{R_{y_0y_0}}(\partial_{u_1}\tilde f_{1}) 
+ (\partial_{u_1} \tilde f_{1})^2\right)(z_0,u_0) = 0.
\]
The three identities in the lemma concerning the second partial derivatives follow immediately.
\end{proof}

\begin{lemma} \label{lem:u1_nolin}
    Let $u_{1,2}(z)$ be the solutions to system (\ref{eq:our_system-2}) and $z_0$ be a the smallest positive $z$ where $u_1(z)$ is singular. 
    If $T \ne 0$ at $(z,u)=(z_0,u_1(z_0))$ then $u_1(z)$ has a square root singularity at $z_0$ and $u_2(z)$ is either analytic or it has a singularity of type $\alpha > 1$.
\end{lemma}

\begin{proof}
    By the computations above we have analytic functions 
$m_0(z,u_1)$,  $m_1(z,u_1)$, and $\tilde  f_1(z,u_1)$ that we can plug into the curve equation.

Now we want to show that for $u_0 =  u_1(z_0)$, it holds that
    \begin{align}
        0 &= u_0^2 -  R_{y_0}(z_0,u_0,\tilde f_1(z_0,u_0), m_1(z_0,u_0), m_0(z_0,u_0)),\label{eq:first}\\
        0 &= 2u_0 - \partial_{u_1}R_{y_0}(z_0,u_0,\tilde f_1(z_0,u_0), m_1(z_0,u_0), m_0(z_0,u_0)),\label{eq:second}\\
        0 &\neq 2-\partial^2_{u_1}R_{y_0}(z_0,u_0,\tilde f_1(z_0,u_0), m_1(z_0,u_0), m_0(z_0,u_0)),\label{eq:third}
    \end{align}
   The square root singularity of $u_1(z)$ will follow  by an application of the Weierstrass preparation theorem and a small argument concerning the valuation of $u$. So, any of the following computations are evaluated at $(z_0,u_0)$. \\
    
    Clearly, equation (\ref{eq:first}) is satisfied, since it is the original curve equation. The second equation (\ref{eq:second}) holds as well, since by Lemma \ref{lem:derivatives}
    \begin{align*}
    &2u_1 - \partial_{u_1}R_{y_0}(z,u_1,\tilde f_1(z,u_1), m_1(z,u_1), m_0(z,u_1))\\
        &\quad =2u_1 - R_{uy_0} - R_{y_0y_0}\partial_{u_1}\tilde f_1 - R_{y_0 y_2} \partial_{u_1}m_0 - R_{y_0 y_1} \partial_{u_1}m_1\\
        &\quad = \,2u_1 - R_{uy_0} + R_{y_0y_0}\frac{2\tilde f_1- R_{uu}}{2u_1-R_{uy_0}}
    \end{align*}
    and the last expression evaluates as $0$ at $(z_0,u_0)$ because $\det B_1 = 0$.
    \\

The third equation (\ref{eq:third}) on the other hand expands to
\begin{align*}
2-\partial^2_{u_1}&R_{y_0}(z,u_1,\tilde f_1(z,u_1), m_1(z,u_1), m_0(z,u_1))\\
=\;2&-R_{uuy_0}-2R_{uy_0y_0}\partial_{u_1} \tilde f_{1}-R_{y_0y_0y_0}(\partial_{u_1} \tilde f_{1})^2 - R_{y_0y_0}\partial^2_{u_1} \tilde f_{1}\\
&-2R_{y_0y_0y_2}(\partial_{u_1} \tilde f_{1})(\partial_{u_1}m_0) - 2R_{y_0y_0y_1}(\partial_{u_1} \tilde f_{1})(\partial_{u_1}m_1) \\
&-2R_{uy_0y_2}(\partial_{u_1}m_0) - 2R_{uy_0y_1}(\partial_{u_1}m_1)\\
&-R_{y_0y_2y_2}(\partial_{u_1}m_0)^2 - R_{y_0y_1y_1}(\partial_{u_1}m_1)^2 - 2R_{y_0y_1y_2}(\partial_{u_1}m_1)(\partial_{u_1}m_0)\\
&- R_{y_0 y_2}\partial^2_{u_1}m_0 - R_{y_0 y_1} \partial^2_{u_1}m_1
\end{align*}
which by Lemma \ref{lem:derivatives} simplifies at $(z_0,u_0)$ to
\begin{align*}
2-\partial^2_{u_1}&R_{y_0}(z_0,u_0,\tilde f_1(z_0,u_0), m_1(z_0,u_0), m_0(z_0,u_0))\\
=
2&-R_{uuy_0}-2R_{uy_0y_0}\partial_{u_1} \tilde f_{1}-R_{y_0y_0y_0}(\partial_{u_1} \tilde f_{1})^2 \\
&- \frac{R_{y_0y_0}}{2u_1-R_{uy_0}}\left(R_{uuu}
+2(R_{uuy_0}-2)\partial_{u_1}\tilde f_{1}+R_{uy_0y_0}(\partial_{u_1}\tilde f_{1})^2\right)
\end{align*}

If we multiply this equation by $\partial_{u_1} \tilde f_{1}(z_0,u_0) = \left(\frac{2u_1-R_{uy_0}}{R_{y_0y_0}}\right)(z_0,u_0)$ (which is non-zero), we see that it is equal to
\[
6\partial_{u_1} \tilde f_{1}-R_{uuu}-3R_{uuy_0}\partial_{u_1}\tilde f_{1}-3R_{uy_0y_0}(\partial_{u_1} \tilde f_{1})^2 
- R_{y_0y_0y_0}(\partial_{u_1} \tilde f_{1})^3 = -T
\]
which is non-zero by assumption.

Now by the Weierstrass preparation theorem, we now know that $u_1(z)$ satisfies a quadratic equation at $z_0$
\[(u_1(z) - u_0)^2 + a_1(z)(u_1(z)-u_0) + a_2(z) = 0\]
where $a_1, a_2$ are analytic functions at $z_0$ and $a_1(z_0) = a_2(z_0) = 0$. By our considerations after Lemma \ref{lem:B_1B_2help} $u_1(z)$ has a singularity of some type $ \alpha \in (0,1)$. If we plug the Puiseux series of $u_1(z)$ at $z_0$ into the equation above, the valuation in $Z = (z-z_0)$ of the first summand in the quadratic equation is $2\alpha$, the valuation of the second is at least $1+\alpha$ because $a_1(z_0)=0$ and the valuation of the last summand denoted by $\beta$ has to be at least 1. Of course the coefficients of the terms in $Z^i$ have to sum up to $0$, so it has to hold that 
\[2\alpha \geq 1 +\alpha \quad \mbox{or} \quad 2\alpha = \beta \geq 1\]
In the first case it would hold $\alpha \geq 1$ which contradicts $\alpha \in (0,1)$. Thus, $2\alpha = \beta$ and since $\beta \in \mathbb{N}$ it follows that $u_1(z)$ has a square root singularity.
\end{proof}

We can determine the singularities of $M_0(z)$ and $M_1(z)$ analogously to the linear case by considering $g(z)$ and $h(z)$ which are both non-zero if the curve equation is not a polynomial in $u^2$ and the split system
\begin{align}
    (g^2+h)f^+ + 2ghf^- + gm_1 + m_0 - R^+(z,g,h, f^+, f^-, m_1, m_0) & = 0, \nonumber \\
(g^2+h)f^- + 2gf^+ + m_1 - R^-(z,g,h, f^+, f^-, m_1, m_0) & = 0 \nonumber  \\
g^2+h - R_{y_0}^+(z,g,h, f^+, f^-, m_1, m_0) &= 0,   \label{eq:split_sys} \\
2g - R_{y_0}^-(z,g,h, f^+, f^-,  m_1, m_0) &= 0,  \nonumber \\
2(gf^++hf^-) + m_1 - R_u^+(z,g,h, f^+, f^-, m_1, m_0) &=0,\nonumber  \\
2 (gf^-+f^+) - R_u^-(z,g,h, f^+, f^-, m_1, m_0) &=0, \nonumber
\end{align}
where $f^+(z) = \left(\diff (z,g\pm h)\right)^+$ and $f^-(z) = \left(\diff (z,g\pm h)\right)^-$. Of course it holds that $f^+(z) = \frac{f_1(z)+f_2(z)}{2}$ and $f^+(z) = \frac{f_1(z)-f_2(z)}{2\sqrt{h}}$.

\begin{lemma} Let $g(z) = z^{d_1} \tilde g(z^d)$ and $h(z)=z^{d_2} \tilde h(z^d)$ be the solutions given by Lemma \ref{Le1b}, then
$\tilde g(z)$ and $\tilde h(z)$ have the (common) radius of convergence $z_0^{1/d}$, where they have a square root singularity.
Furthermore there exists $\delta> 0$ such that $\tilde g(z)$ and $\tilde h(z)$ can be analytically continued to
the region $\mid z \mid < z_0^{1/d} + \delta$, $\mid z-z_0^{1/d} \mid > 2\delta$.
\end{lemma} 

\begin{proof}
Since $u_1(z)$ has a square root singularity and $u_2(z)$ a singularity of type $\alpha> 1$, the functions $g(z)$ and $h(z)$ also have a common square root singularity at $z=z_0$.
Once again, as a consequence the functions $\tilde g(z)$ and $\tilde h(z)$ have a common square root singularity 
at $z = z_0^{1/d}$ which in turn generates the square root singularities of $g(z)$ and $h(z)$ at
$z = z_0 e^{2\pi i \ell/d}$, $1\le \ell < d$.\\

Now, we have a look at the Jacobian of the split system (\ref{eq:split_sys}) and show that its determinant only vanishes for $z = z_0 e^{2\pi i \ell/d}$, $1\le \ell < d$.

Observe that
\begin{align*} (R^+)_g &=  \frac{R_u+\overline{R_u}}{2} = (R_u)^+ \\
(R^-)_g &=  \frac{R_u-\overline{R_u}}{2\sqrt{h}} = (R_u)^- \\
(R^+)_h &= \frac{R_u-\overline{R_u}}{4\sqrt{h}}+\frac{R_{y_0}-\overline{R_{y_0}}}{4\sqrt{h}}f^- \\
&= \frac{1}{2} ((R_u)^-+R_{y_0}^-f^-) = 2gf^- + f^+\\
(R^-)_h &= \frac{R_u+\overline{R_u}}{4h}+ \frac{R_{y_0}+\overline{R_{y_0}}}{4h}f^- - \frac{R-\overline{R}}{4h\sqrt{h}} \\
&= \frac{1}{2h}\left((R_u)^++(R_{y_0})^+f^--R^-\right) = f^-\\
(R^+)_{f^+} &=  \frac{R_{y_0}+\overline{R_{y_0}}}{2} = (R_{y_0})^+\\
(R^-)_{f^+} &=  \frac{R_{y_0}-\overline{R_{y_0}}}{2\sqrt{h}} = (R_{y_0})^-\\
(R^+)_{f^-} &=  \frac{R_{y_0}-\overline{R_{y_0}}}{2}\sqrt{h} = (R_{y_0})^-h\\
(R^-)_{f^-} &=  \frac{R_{y_0}+\overline{R_{y_0}}}{2} = (R_{y_0})^+
\end{align*}
Thus, the Jacobian of (\ref{eq:split_sys})  with respect to $ m_1, m_0, g, h , f^+, f^-$ equals
\[\left(\begin{matrix}
    A' & 0 \\
    C' & B'
\end{matrix}\right)\]
for matrices $A' \in \mathbb{C}^{2\times 2}, B' \in \mathbb{C}^{4\times 4}$:
\[A' = \left(\begin{matrix}
    1-R_{y_2}^+ & g-R_{y_1}^+ \\
    -R_{y_2}^- & 1-R_{y_1}^-
\end{matrix}\right)\]
\[B' = \left(\begin{matrix}
    2g-R_{uy_0}^+ & 1-\frac{1}{2}(R^-_{uy_0}+R_{y_0y_0}^-f^-) & -R^+_{y_0y_0} & -hR^-_{y_0y_0}\\
    2-R^-_{uy_0} & -\frac{1}{2h}(R^+_{uy_0}+R_{y_0y_0}^+f^--R^-_{y_0}) & -R^-_{y_0y_0} & -R^+_{y_0y_0}\\
    2f^+-R^+_{uu} & 2f^--\frac{1}{2}(R^-_{uu}+R_{uy_0}^-f^-) & 2g-R^+_{uy_0} & 2h-hR^-_{uy_0}\\
    2f^--R^-_{uu} & -\frac{1}{2h}(R_{uu}^++R_{uy_0}^+f^--R_u^-) & 2-R_{uy_0}^- & 2g-R_{uy_0}^+,
\end{matrix}  \right)\]
and a matrix $C' \in \mathbb{C}^{4\times 2}$.
The zeroes arise again because the equations
\[\left(\begin{matrix}
   2gf^++2hf^-+m_1-R_u^+  & f^++2gf^--R^+_h & g^2+h-R^+_{y_0} & 2gh-hR_{y_0}^-\\
   2gf^-+2f^+-R_u^- & f^--R^-_h & 2g-R^-_{y_0} & g^2+h-R_{y_0}^+
\end{matrix}\right) = 0.\]

The determinant thus factors into the determinants of the matrices $A'$ and $B'$.
First, we observe that  
\begin{align*}
    \det \left(\begin{matrix}
    1-R_{y_2}^+ & g-R_{y_1}^+ \\
    -R_{y_2}^- & 1-R_{y_1}^-
    \end{matrix}\right) 
    &= \det \left(\begin{matrix}
    1-R_{y_2}^+ & -(R_{y_1}^+-gR_{y_2}^+) \\
    -R_{y_2}^- & 1-(R_{y_1}^--gR_{y_2}^-)
    \end{matrix}\right) \\
    &= \det \left(\begin{matrix}
    1-R_{y_2}^+ & -(zQ_{\alpha_0}^-h + zQ_{\alpha_1}^+) \\
    -R_{y_2}^- & 1-(zQ_{\alpha_0}^+ +zQ_{\alpha_1}^-)
    \end{matrix}\right) \\
    &= (1-R_{y_2}^+)\left(1-(zQ_{\alpha_0}^+ + zQ_{\alpha_1}^-)-\frac{(hzQ_{\alpha_0}^-+zQ_{\alpha_1}^+)R_{y_2}^-}{1-R_{y_2}^+} \right)
\end{align*}
where
\[D_2(z,g(z),h(z)) = (zQ_{\alpha_0}^++zQ_{\alpha_1}^-)+\frac{(hzQ_{\alpha_0}^-+zQ_{\alpha_1}^+)R_{y_2}^-}{1-R_{y_2}^+}\]
is a power series with non-negative coefficients. Now, given that $h(z) > 0$ for $0<z\leq z_0$, we show that $\mid D_2(z,g(z),h(z))\mid < 1$
for all $\mid z \mid < z_0$ by using the relation
\begin{align*}\det A' &= \frac{1}{2\sqrt{h}}\det \left(\begin{matrix}
    1-R_{y_2}^+-R_{y_2}^-\sqrt{h} & g-R_{y_1}^+ +(1-R_{y_1}^-)\sqrt{h}\\
    -2R_{y_2}^-\sqrt{h} & -2(1-R_{y_1}^-)\sqrt{h}
\end{matrix}\right) \\
&= \frac{-1}{2\sqrt{h}}\det \left(\begin{matrix}
    1-R_{y_2} & u_1-R_{y_1}\\
    1-\overline{R_{y_2}} & u_2-\overline{R_{y_1}}
\end{matrix}\right) = \frac{-1}{2\sqrt{h}}\det A. \end{align*}
and the fact that $\det A \ne 0$ at $z_0$.\\

Similarly, we can relate the determinant of $B'$ to the determinants of $B_1$ and $B_2$.
First we multiply the second and fourth line and the second column of $B'$ by $\sqrt{h}$. Then we divide the last column by $\sqrt{h}$ and add the last column of the obtained matrix multiplied by $f^-$ to the second column. Finally, we multiply the second column by $2$ and obtain

\begin{align*}
    \det B' &= \frac{1}{2h}\det \left(\begin{matrix}
    2g-R_{uy_0}^+ & (2-R^-_{uy_0})\sqrt{h} & -R^+_{y_0y_0} & -R^-_{y_0y_0}\sqrt{h}\\
    (2-R^-_{uy_0})\sqrt{h} & -(R^+_{uy_0}-R^-_{y_0}) & -R^-_{y_0y_0}\sqrt{h} & -R^+_{y_0y_0}\\
    2f^+-R^+_{uu} & (2f^--R^-_{uu})\sqrt{h} & 2g-R^+_{uy_0} & (2-R^-_{uy_0})\sqrt{h}\\
    (2f^--R^-_{uu})\sqrt{h} & -(R_{uu}^++2gf^--R_u^-) & (2-R_{uy_0}^-)\sqrt{h} & 2g-R_{uy_0}^+
\end{matrix}  \right) \\
&= \frac{1}{2h}\det \left(\begin{matrix}
    2g-R_{uy_0}^+ & (2-R^-_{uy_0})\sqrt{h} & -R^+_{y_0y_0} & -R^-_{y_0y_0}\sqrt{h}\\
    (2-R^-_{uy_0})\sqrt{h} & 2g-R^+_{uy_0} & -R^-_{y_0y_0}\sqrt{h} & -R^+_{y_0y_0}\\
    2f^+-R^+_{uu} & (2f^--R^-_{uu})\sqrt{h} & 2g-R^+_{uy_0} & (2-R^-_{uy_0})\sqrt{h}\\
    (2f^--R^-_{uu})\sqrt{h} & 2f^+-R_{uu}^+ & (2-R_{uy_0}^-)\sqrt{h} & 2g-R_{uy_0}^+
\end{matrix}  \right)\\
&= \frac{1}{8h}\det \left(\begin{matrix}
    2u_1-R_{uy_0} & 2u_1-R_{uy_0} & -R_{y_0y_0} & -R_{y_0y_0}\\
    2u_2-\overline{R_{uy_0}} & -(2u_2-\overline{R_{uy_0}}) & -\overline{R_{y_0y_0}} & \overline{R_{y_0y_0}}\\
    2\diff_1-R_{uu} & 2\diff_1-R_{uu} & 2u_1-R_{uy_0} & 2u_1-R_{uy_0}\\
    2\diff_2-\overline{R_{uu}} & -(2\diff_2-\overline{R_{uu}}) & 2u_2-\overline{R_{uy_0}} & -(2u_2-\overline{R_{uy_0}})
\end{matrix}  \right)\\
&= \frac{1}{8h}\det \left(\begin{matrix}
    0 & 2u_1-R_{uy_0} & 0 & -R_{y_0y_0}\\
    2u_2-\overline{R_{uy_0}} & 0 & -\overline{R_{y_0y_0}} & 0\\
    0 & 2\diff_1-R_{uu} & 0 & 2u_1-R_{uy_0}\\
    2\diff_2-\overline{R_{uu}} & 0 & 2u_2-\overline{R_{uy_0}} & 0
\end{matrix}  \right)\\
& = \frac{1}{8h} \det B_1 \det B_2
\end{align*}
Hence, the determinant of the Jacobian equals the determinant of the original system up to a factor of $16h^{3/2}$. Again, we want the zeroes of this Jacobian be very local such that we can analytically continue $g(z)$ and $h(z)$ in a $\Delta$-region around them.\\
To show this we use another (not very surprising) relation arising from the system $g^2+h = C^+(z,g,h), 2g =C^-(z,g,h)$. By Lemma \ref{lem:B_1B_2help}, it holds that

\begin{align*} 
\det \left(\begin{matrix}
    1-C_h^+ & 2g - C_g^+\\
    -C_h^- & 2-C_g^-
\end{matrix}\right) &= \frac{1}{2h}\det \left(\begin{matrix}
    (2-C_u^-)\sqrt{h} & 2g - C_u^+\\
    2g-C_u^+ & (2-C_u^-)\sqrt{h}
\end{matrix}\right) \\
&=   -\frac{(2u -C_u)(\overline{2u-C_u})}{2h} = \frac{\sqrt{\det B_1 \det B_2}}{2h}\\
&= 2(1-C_h^+)\left(1-C_g^-/2 - \frac{(C_g^+ - C^-) C_h^-/2}{1- C_h^+}\right).
\end{align*}
Again we consider the power series with non-negative coefficients
\[D(z,g,h) = C_g^-/2 + \frac{(C_g^+ - C^-) C_h^-/2}{1- C_h^+}\]
which has to equal 1 at $z_0$ by continuity. Analogously to the linear case, it holds for
$\mid z \mid = z_0$ but $z\ne z_0 e^{2\pi i \ell/d}$, $0\le \ell < d$ that
\[
\mid g(z) \mid < g(z_0) \quad \mbox{and}\quad \mid h(z) \mid < h(z_0)
\]
and we can therefore conclude that
\[
\mid D(z,g(z),h(z)) \mid  \le D(\mid z \mid, \mid g(z_0) \mid, \mid h(z_0) \mid) 
< D(z_0,g(z_0),h(z_0)) = 1.
\]
Hence the factor $\det B_1 \det B_2$ of the determinant of the Jacobian of the split system is non-zero in a neighborhood of $z_0$ excluding $z_0$.

Thus,
$g(z)$ and $h(z)$ can be analytically continued in a $\Delta$-region.
\end{proof}

\begin{lemma}\label{le:Le5}
Suppose that the assumptions of Theorem~\ref{Th2} are satisfied.  
Then $M_0(z)$ and $M_1(z)$ have a 3/2-singularity at $z_0$.
Furthermore there (uniquely) exists an integer $d\ge 1$ 
such that $M_0(z)$ and $M_1(z)$ can be analytically continued to the region
\begin{equation}\label{eqanalyticcontregion}
\mid z \mid < z_0 + \delta, \quad  \mid z - z_0 e^{2\pi i \ell/d} \mid > 2 \delta \quad (0\le \ell < d)
\end{equation}
for any sufficiently small $\delta > 0$ and both functions have a 3/2-singularity at $z_0$ and at the points $z_0 e^{2\pi i \ell/d}$, $1\le \ell < d$.
\end{lemma}

\begin{proof}
If we only consider the first, second, fifth and sixth line of the split system (\ref{eq:split_sys}), the Jacobian with respect to $m_0, m_1, f^+, f^-$ is invertible at $z_0$ and we obtain analytic functions $\hat m_0(z,g,h), \hat m_1(z,g,h), \hat f^+(z,g,h), \hat f^-(z,g,h)$ for which it holds that
\begin{align*}M_0(z) = \hat m_0(z,g(z),h(z)), \quad M_1(z) = \hat m_1(z,g(z),h(z)).
\end{align*}
Again, we will compute a Taylor expansion to confirm the $3/2$-singularities arising from the singularities of $g$ and $h$.
\[\partial_g \left(\begin{matrix} m_0 \\ m_1 \end{matrix} \right) = (A')^{-1} \left(\begin{matrix}
    2gf^++2hf^-+m_1-R^+_u \\
    2gf^-+2f^+-R^-_u
\end{matrix}\right) = 0.\]
Similarly,
\[\partial_h \left(\begin{matrix} m_0 \\ m_1 \end{matrix} \right) = (A')^{-1} \left(\begin{matrix}
    f^++2gf^--(2gf^-+f^+) \\
    f^--f^-
\end{matrix}\right) = 0.\]
That means 
\begin{align*}M_0(z) = \hat m_0(z,g(z_0),h(z_0)) + O(z-z_0), \quad M_1(z) = \hat m_1(z,g(z_0),h(z_0))+ O(z-z_0).
\end{align*}
and $M_0(z)$ and $M_1(z)$ have at most 3/2 singularities at $z_0e^{2\pi i/\ell}$, $0\le \ell < d$.
However, analogously to the linear case, we can express $M_1'(z)$ as a power series with positive coefficients in $g,h$ along the solutions
\[M_1'(z) = \left(\frac{R_z^-+R_{y_2}M_0'(z)}{1-R_{y_1}}\right)(z,g(z),h(z),f^+(z),f^-(z),M_1(z),M_0(z))\]
which confirms square root singularities at $z_0e^{2\pi i/\ell}$, $0\le \ell < d$ and therefore $3/2$-singularities for $M_1(z)$. 
Furthermore, by applying Lemma~\ref{Le3} again it follows that 
$M_1(z)$ can be analytically continued to a region of the form (\ref{eqanalyticcontregion})
for any sufficiently small $\delta > 0$.\\

Finally, since the system is strongly connected the same holds for $M_0(z)$.
\end{proof}

\begin{lemma}\label{le:Le6}
Suppose that the assumptions of Theorem~\ref{Th2} are satisfied and let $d\ge 1$ be
the (unique) integer given in Lemma~\ref{le:Le5}.  
Then there exists a non-empty set $J \subseteq \{0,1,\ldots,d-1\}$ 
of residue classes modulo $d$ and constants $c_j> 0$ such that for $j\in J$
\begin{equation}\label{eqTh1.2}
M_n = [z^n]\,M_0(z) = c_j n^{-3/2} z_0^{-n}\left( 1 + O\left( \frac 1n \right) \right), \qquad (n\equiv j \bmod b,\ n\to \infty).
\end{equation}
Furthermore, if $n\equiv j \bmod d$ with $j\not \in J$ then we either have $M_n = 0$ 
or
\[
M_n = O\left( n^{-5/2} z_0^{-n} \right).
\]
\end{lemma}

\begin{proof}
Let $J_1$ denote the set of $\ell \in \{0,\ldots, d-1\}$ for which
$z= z_0 e^{2\pi i \ell/d}$ is a 3/2-singularity of $M_0(z)$. $J_1$ is non-empty since $0 \in J_1$.
Suppose that the local expansions around $z= z_0 e^{2\pi i \ell/d}$ (for $\ell \in J_1$) are given by
\[
M_0(z) = g_\ell(z) - h_\ell(z) \left(1- \frac{z}{z_0}  e^{-2\pi i \ell/d} \right)^{3/2}
\]
for certain functions $g_\ell,h_\ell$ that are analytic at $z= z_0 e^{2\pi i \ell/d}$.
The by using a standard singularity analysis (see \cite{FO}) together with the properties stated in Lemma~\ref{Le5} it follows that the $n$-th coefficient $M_n = [z^n]\,M_0(z)$ is asymptotically given by 
\[
M_n =  \frac 1{2\sqrt \pi } \sum_{\ell \in J_1} b_\ell  e^{-2\pi i \ell n/d} n^{-5/2} z_0^{-n} 
+ O\left( n^{-7/2} z_0^{-n} \right)
\]
where $b_\ell = h_\ell(z_0 e^{2\pi i \ell/d})$. Hence, if $n \equiv j \bmod d$ then we have
\[
M_n = c_j n^{-5/2} z_0^{-n} + O\left( n^{-7/2} z_0^{-n} \right),
\]
where
\[
c_j = \frac 1{2\sqrt \pi } \sum_{\ell \in J_1} b_\ell  e^{-2\pi i j \ell/d}.
\]
Since the coefficients $M_n$ are non-negative it follows that the numbers $c_j$ are non-negative, too.
The set $J = \{ j\in \{0,\ldots d-1\} : c_j > 0\}$ satisfies the desired result.
\end{proof}

\section{Proof of Theorem~\ref{Th3}}

Suppose first that we fix $w > 0$. Then the equation (\ref{eqMzuw}) satisfies either the assumptions of
Theorem~\ref{Th1} or Theorem~\ref{Th2}. Actually $w$ can be considered just as a {\it scaling} of the
equation but does not change the structure of the solution. In particular the singularity structure
of $M(z,0,w)$ does not depend on $w$, only the location of the dominant singularity $z_0(w)$ and 
constant evaluations depend on $w$. Furthermore if $w$ is sufficiently close to $1$ then the condition
$T\ne 0$ (in Theorem~\ref{Th2}) will still hold. Hence, there exist $d\ge 1$ and a non-zero set 
$J$ of residue classed modulo $d$ (both of them indepdendent of $w$) such that for $j\in J$
and as $n\to\infty$  ($n\equiv j \bmod b$)
\[
[z^n]\, M(z,0,w) 
= c_j(w) n^{-3/2} z_0(w)^{-n}\left( 1 + O\left( \frac 1n \right) \right)
\]
or
\[
[z^n]\, M(z,0,w) = c_j(w) n^{-5/2} z_0(w)^{-n}\left( 1 + O\left( \frac 1n \right) \right).
\]
Actually these asymptotic expansions hold uniformly in $w$ if $w$ varies in a compact interval
of the positive reals (compare, for example,  with \cite[Lemma~2.18]{D}).
Hence it follows that the probability generating function of $X_n$ ($n\equiv j \bmod b$) is asympotically given by
\[
\mathbb{E}\, w^{X_n}  = \frac{ [z^n]\, M(z,0,w) } { [z^n]\, M(z,0,1) } 
= \frac{c_j(w)  }{c_j(1)} 
\left(  \frac {z_0(1)}  {z_0(w)}   \right)^n 
\left( 1 + O\left( \frac 1n \right) \right), 
\]
Hence, by a direct application of Hwang's quasi-power theorem \cite[Theorem~2.22]{D} it follows 
that $X_n$ ($n\equiv j \bmod b$) satisfies a central limit theorem as proposed.

\section{Examples}\label{sec6}

In this section, we will illustrate our generic computations in the proof of Theorems \ref{Th1} and \ref{Th2} on the examples given in Section \ref{sec:generic}.

\begin{example} [Example~\ref{ex3} continued]
For one-dimensional non-negative lattice paths where we allow steps of the form $\pm 1$ and $\pm 2$ we obtained the functional equation
\begin{equation}\label{eqEzv2}
E(z,u) = 1 + z(u+u^2) E(z,u) + z \frac{E(z,u) - E(z,0)}u + z \frac{E(z,u) - E(z,0) - u E_v(u,0) }{u^2}.
\end{equation}
We know that the curve equation
\[u^{2} = z(1+u)u^{3}+zu+z.\]
has two solutions $u_1(z)$, $u_2(z)$ with $u_1(0) = u_2(0) = 0$ and which are real valued along the positive reals. In order to determine $z_0>0$ where $u_1(z)$ is singular, we consider the partial derivative with respect to $u$
\[2u = 3zu^2+4zu^3+z.\]
The common zeros $(z_{0},u_{0})$ of this equation and the curve equation are
\[\left \{(0,0), \left(\frac{1}{4},1\right), \left(-\frac{4}{9}, -\frac{1-\sqrt{15}i}{4} \right),\left(-\frac{4}{9}, -\frac{1+\sqrt{15}i}{4} \right)\right \}\]
Hence, it follows that $z_0 = \frac 14$ and $u_1(z_0) = 1$. Furthermore, the local expansion of $u_1(z)$ at $z=z_0$ is 
given by
\[
u_1(z) = 1 - \sqrt{8} \sqrt{1-4z} + \cdots.
\]

If we want to compute the local expansion of $M_0(z)$ and $M_1(z)$, we also need the value of $u_2(z_0)$. Of course the equation
\[ u_2(z_0)^{2} = \frac 14(1+u_2(z_0))u_2(z_0)^{3}+\frac 14u_2(z_0)+\frac 14 \]
has to be satisfied and admits the following possible values for $u_2(z)$
\[\left\{1, -\frac 32 \pm \frac{\sqrt{5}}{2}\right\}.\]
Since $u_1(z_{0}) = 1$  and it has to hold that $\mid u_2(z_0)\mid < u_1(z_0)$, the correct value is 
\[
u_2\left(z_{0} \right) = \frac{\sqrt 5 - 3}2.
\]
Alternatively, we can also consider the system of equations for $g(z)$ and $h(z)$:
\begin{align*}
g^2 + h &=  z(g^4 + h^2 + g^3 + 3(2g^2 + g)h + g + 1), \\
2g  &= z((4g + 1)h + (4g^3 + 3g^2 + 1)).
\end{align*}
At $z_0 = \frac 14$ there are only finitely many solutions, however, the only ones with $g_0 + \sqrt{h_0} = 1$
are $g_0 = (\sqrt 5 -1)/4$ and $h_0 = (15 - 5\sqrt 5)/8$. Consequently we have
\[
u_2(z_0) = g_0 - \sqrt{h_0} = \frac{\sqrt 5 - 3}2.
\]
Finally we use the linear system (\ref{eqM01-1})--(\ref{eqM01-2}) to obtain the local expansion for
$M_0(z)$ and $M_1(z)$:
\begin{align*}
M_0(z) &=  (6- 2 \sqrt 5) - \frac{101\sqrt 2 - 45 \sqrt{10} }{19} \sqrt{1-4z} + \cdots  ,\\
M_1(z) &= (4 \sqrt 5 - 8) -  \frac{28\sqrt {10} - 62 \sqrt{2} }{19} \sqrt{1-4z} + \cdots.
\end{align*}
Since there are no periodicities this implies that
\[
[z^{n}]E(z,0) \sim  \frac{101\sqrt 2 - 45 \sqrt{10} }{38 \sqrt{\pi}} \,  n^{-3/2}\,4^{n}
\]
\end{example}

\begin{example}  [Example~\ref{ex4} continued]
The functional equation for $3$-Constellations can be transformed to the equation
\begin{align*}
\tilde C(z,v) &= z(v+1)(\tilde C(z,v)+1)^3 + z(v+1)(2 \tilde C(z,v) + \tilde C(z,0)+3) \frac{\tilde C(z,v)- \tilde C(z,v)}{v}  \\
&+ z(v+1)\frac{\tilde C(z,v)- \tilde C(z,0)- v\tilde C_v(z,0)}{v^2}
\end{align*}
by substituting $v=u-1$ and defining $\tilde C(z,v) =C(z,u)-1$. Note that of course $\tilde C(z,0) = C(z,1)-1$ and $\tilde C_v(z,0) = C_v(z,1) = C_u(z,1)$. The equations for the unknown functions $v_{1}(z), v_2(z), d_{1} = \diff^{(2)}\tilde C(z,v_{1}(z)), d_{2} = \diff^{(2)}\tilde C(z,v_{2}(z))$ and $m_{1}(z) = \tilde C_{v}(z,0)$, $m_{0} = \tilde C(z,0)$ are thus given by
\begin{align*}
v_{i}^2d_{i} +v_{i}m_{1}+ m_{0} =  &z(v_{i}+1)(v_{i}^2d_{i} +v_{i}m_{1}+ m_{0}+1)^3 \\
&+ z(v_{i}+1)((2v_{i}^2d_{i} +2v_{i}m_{1}+ 3m_{0} + 3)(v_{i}d_{i} +m_{1}) + d_{i})\\
v_{i}^{2} = &z(v_{i}+1)(3v_{i}^2(v_{i}^2d_{i} +v_{i}m_{1}+ m_{0}+1)^2  + 2v_{i}^2(v_{i}d_{i} +m_{1})) \\
&+ z(v_{i}+1)((2v_{i}^2d_{i} +2v_{i}m_{1}+ 3m_{0} + 3)v_{i} + 1),\\
2v_{i}d_{i}+m_{1} = &z(v_{i}^2d_{i} +v_{i}m_{1}+ m_{0}+1)^3\\
&+z(2v_{i}^2d_{i} +2v_{i}m_{1}+ 3m_{0} + 3)(v_{i}d_{i} +m_{1}) + zd_{i} \\
&+3z(v_{i}+1)(v_{i}^2d_{i} +v_{i}m_{1}+ m_{0} + 1)^2(2v_{i}d_{i}+m_{1} ) \\
&+ 2z(v_{i}+1)(2v_{i}d_{i}+m_{1} )(v_{i}d_{i}+m_{1} ) \\
&+ z(v_{i}+1)(2v_{i}^2d_{i} +2v_{i}m_{1}+ 3m_{0}  + 3)d_{i}. \end{align*}

Numerical computations show that the smallest positive $z_{0}$ where the Jacobian of this system is invertible equals
$z_{0} \approx 0.0494$.  Indeed, the exact value for the singularity is $4/81 = 0,04938...$. The other variables take the 
approximate values
\[u_{1} \approx  0.6867, \quad u_{2} \approx  -0.1562, \quad d_{1} \approx  0.1070, \quad d_{2} \approx  0.0433, \quad m_{1} \approx 0.1134,\quad m_{0} \approx 0.0833.\]
Note that all computations can be worked out although the scheme of Theorem~\ref{Th2} is not strictly satisfied. 
The values of the determinants at $z_0$ equal
\[\det A \approx  -0.2588, \quad \det B_{1} \approx  0, \quad \det B_{2} \approx  0.1828\]
such that we can apply Lemma \ref{lem:B_1B_2help}.
The neccessary condition of Theorem~\ref{Th2} is also satisfied, since the value of the expression equals $T \approx 2.7209$.
Finally we get the asymptotics
\[[z_{n}] C(z,0) \sim c\,n^{-5/2}\left(\frac{81}{4}\right)^{n}\]
for $c \approx 0.0731$.
\end{example}

\bibliography{universalcatalytic}{}
\bibliographystyle{plain}

\end{document}